\documentclass[a4paper,reqno,12pt]{amsart}
\usepackage[all]{xy}           %For commutative diagrams
\usepackage{amssymb}           %For double arrows
\usepackage[utf8]{inputenc}
\usepackage{hyperref}
\usepackage{eucal}
\usepackage{epsfig}
\usepackage{pst-grad} % For gradients
\usepackage{pst-plot} % For axes
\numberwithin{equation}{section}

\newtheorem{definition}{Definition}[section]
\newtheorem{lemma}[definition]{Lemma}
\newtheorem{theorem}[definition]{Theorem}
\newtheorem{proposition}[definition]{Proposition}

\newtheorem{remarkth}[definition]{Remark}
\newtheorem{example}[definition]{Example}

\renewcommand{\emph}[1]{{\bfseries\itshape{#1}}}

\makeatletter
\newcommand\prol{\@ifstar{\@proldf}{\@prolpf}}  %% if * dual else primal
\def\@prolpf{\@ifnextchar[{\@prolpf@wrt}{\@prolpf@}}
\def\@prolpf@wrt[#1]#2{\@ifnextchar[{\@prolpf@wrt@at{#1}{#2}}{\@prolpf@wrt@{#1}{#2}}}
\def\@prolpf@wrt@at#1#2[#3]{\prolsymbol^{#1}_{#3}#2}
\def\@prolpf@wrt@#1#2{\prolsymbol^{#1}#2}
\def\@prolpf@#1{\@ifnextchar[{\@prolpf@at{#1}}{\@prolpf@@{#1}}}
\def\@prolpf@at#1[#2]{\prolsymbol_{#2}#1}
\def\@prolpf@@#1{\prolsymbol#1}
\def\@proldf{\@ifnextchar[{\@proldf@wrt}{\@proldf@}}
\def\@proldf@wrt[#1]#2{\@ifnextchar[{\@proldf@wrt@at{#1}{#2}}{\@proldf@wrt@{#1}{#2}}}
\def\@proldf@wrt@at#1#2[#3]{\prolsymbol^{*#1}_{#3}#2}
\def\@proldf@wrt@#1#2{\prolsymbol^{*#1}#2}
\def\@proldf@#1{\@ifnextchar[{\@proldf@at{#1}}{\@proldf@@{#1}}}
\def\@proldf@at#1[#2]{\prolsymbol^*_{#2}#1}
\def\@proldf@@#1{\prolsymbol^*#1}
\def\prolsymbol{\mathcal{T}}
\makeatother

\begin{document}

\title[Groupoids and algebroids in uniformity and homogeneity]{Lie groupoids and algebroids applied to the study of uniformity and homogeneity of material bodies}

\author[V. M. Jim\'enez]{V\'ictor Manuel Jim\'enez}
\address{V\'ictor Manuel Jim\'enez:
Instituto de Ciencias Matem\'aticas (CSIC-UAM-UC3M-UCM),
c$\backslash$ Nicol\'as Cabrera, 13-15, Campus Cantoblanco, UAM
28049 Madrid, Spain} \email{victor.jimenez@icmat.es}

\author[M. de Le\'on]{Manuel de Le\'on}
\address{Manuel de Le\'on: Instituto de Ciencias Matem\'aticas (CSIC-UAM-UC3M-UCM),
c$\backslash$ Nicol\'as Cabrera, 13-15, Campus Cantoblanco, UAM
28049 Madrid, Spain.
\newline
and
\newline
Real Academia de Ciencias Exactas, F{\'i}sicas y Naturales, c$\backslash$ de Valverde, 22, 28004 Madrid, Spain.
} \email{mdeleon@icmat.es}

\author[M. Epstein]{Marcelo Epstein}
\address{Marcelo Epstein:
Department of Mechanical Engineering. University of Calgary. 2500 University Drive NW, Calgary, Alberta, Canada, T2N IN4} \email{epstein@enme.ucalgary.ca}

\keywords{Lie algebroid, Lie groupoid, homogeneity, uniformity, material groupoid, material algebroid, $G-$structure, derivation}

 \subjclass[2000]{}

\begin{abstract}
A Lie groupoid, called \textit{material Lie groupoid}, is associated in a natural way to any elastic material. The corresponding Lie algebroid, called \textit{material algebroid}, is used to characterize the uniformity and the homogeneity properties of the material. The relation to previous results in terms of $G-$structures is discussed in detail. An illustrative example is presented as an application of the theory.
\end{abstract}

\maketitle

\tableofcontents

\section{Introduction}

In Continuum Mechanics a simple material body $\mathcal{B}$ is represented by a three-dimensional differentiable manifold which can be covered with just one chart (see for example \cite{JEMARS,GAMAU2,CTRUE,CCWAN}). Given a simple material body $\mathcal{B}$ we identify an embbeding $\phi_{0} : \mathcal{B} \rightarrow \mathbb{R}^{3}$ as a reference configuration. A Lie groupoid, called \textit{material groupoid}, can be naturally associated to any material body (see for example \cite{MEPMDLSEG}, or even \cite{MEPMDL,COSVME} for Cosserat media). In particular, given two different points $X$ and $Y$ of the material body $\mathcal{B}$, a material isomorphism is a linear isomorphism $P_{XY}: T_{X}\mathcal{B} \rightarrow T_{Y} \mathcal{B}$ such that the mechanical response at $X$ and $Y$ is the same; more precisely, if $W$ is the response functional depending on the deformation gradient $F$ at any point $X$ of $\mathcal{B}$, then,
$$ W \left( F P , X \right) = W \left( F , Y \right),$$
for any $F$. The collection of all material isomorphisms for all pairs of points of $\mathcal{B}$ is just the material groupoid, which will be denoted by $\Omega \left( \mathcal{B} \right)$.\\
The theory of continuous distributions of defects based on the notion of material isomorphism has been developed by W. Noll \cite{WNOLL} (see also \cite{CTRUETOU,CTRUE,CCWAN,CCWANSEG}). A structurally based theory had been originally conceived by K. Kondo \cite{KKON}, D. A. Bilby \cite{DABIL}, E. Kröner \cite{EKRONFI,EKRON}, J.D. Eshelby \cite{JDES} and others, as the result of a limiting process starting from a defective crystalline structure (see also the books by R. Lardner \cite{RLAND} and F.R. Navarro \cite{FRNAV}).\\
 In Noll's terminology, a body is {\it uniform} if all its points are mutually materially isomorphic. In generalizing the work of Noll, the notion of material groupoid encodes all the information pertaining to the local and distant symmetries of a material body. In particular, material uniformity corresponds to the {\it transitivity} of the material groupoid. Since transitive groupoids can be regarded as $G-$structures, previous treatment in terms of $G-$structures were necessarily limited to uniform bodies. This severe limitation is removed when using an approach based on the theory of groupoids, thus making the theory applicable to non-uniform bodies, such as functionally graded materials.\\
Another crucial point about this material groupoid is the availability of the associated Lie algebroid, $A \Omega \left( \mathcal{B} \right)$, which is the infinitesimal version of $\Omega \left( \mathcal{B} \right)$.\\
A uniform body is said to be {\it locally homogeneous} if for each body point there exists a configuration whereby the Euclidean translations in an open neighborhood of the point are material isomorphisms. It is a remarkable fact that this homogeneity can be characterized through the properties of $A \Omega \left( \mathcal{B} \right)$. This is indeed accomplished, and related with the earlier approach developed in \cite{MELZA} (see also \cite{FBLOO}) in the framework of $G-$structures.\\
Our paper is divided in three parts: The first part (Sections 1 and 2) is a brief introduction to the fundamental concepts which we will need in the paper. From the mathematical perspective, these concepts are: \textit{Groupoids} and \textit{Lie algebroids}. The fundamental physical notions are \textit{uniformity} and \textit{homogeneity} of simple bodies.\\
Section 3 is devoted to a purely mathematic development. We introduce the notion of \textit{integrability} of Lie subgroupoids of the $1-$jet Lie groupoid on a manifold $M$ ($\Pi^{1} \left( M , M \right)$). We deal with this notion and we relate it with the corresponding one in the $G-$structures and Lie algebroids contexts.\\
In Section 4 we use the results of Section 3 to characterize the uniformity and the homogeneity of simple material bodies. Finally, we use these results to deal with a class of examples, namely, \textit{simple liquid crystals} presented by Coleman \cite{COLE} and Wang \cite{CCWANTHIRD}. These are some of the most common materials that, in the continuous limit, can be, and have been, modeled as \textit{Cosserat media}. An example of the study of dislocations in these media can be found in \cite{MEPMDLTHIRD}. It is possible also to develop a theory of liquid crystals within the context of simple media at the price of losing some of their distinctive physical properties \cite{COLE}. One of our purposes in this example is to demonstrate how the material algebroid can be constructed directly from the constitutive equation in the case of simple materials (a model that may be apt for the description of certain type A smectic materials).
\section{Groupoids and Lie algebroids}
\subsection*{Groupoids}
In this first section we provide a brief introduction to the notions of \textit{Lie groupoid} and \textit{Lie algebroid}. For details we mainly refer to \cite{KMG}. There are also good introductions to groupoids in \cite{EPSBOOK} and \cite{WEINSGROUP}. Another recommendable book as an introduction to these two topics is the book \cite{JNM} (in Spanish).\\
Roughly speaking, the notion of \textit{groupoid} is similar to the notion of group, although the composition is not totally but only partially defined. More specifically, a groupoid is given by two sets, $M$ (\textit{base}) and $\Gamma$ (\textit{total space}), provided with the maps $\alpha,\beta : \Gamma \rightarrow M$ (\textit{source map} and \textit{target map} respectively), $\epsilon: M \rightarrow \Gamma$ (\textit{identities map}), $i: \Gamma \rightarrow \Gamma$ (\textit{inversion map}) and   
$\cdot : \Gamma_{\left(2\right)} \rightarrow \Gamma$ (\textit{composition law}) where,
$$\Gamma_{\left(2\right)} := \{ \left(g,h\right) \in \Gamma \times \Gamma \ : \ \alpha\left(g\right)=\beta\left(h\right) \}.$$
The composition law should satisfy the associative property and the identities and the inversions satisfy that for all $ g \in \Gamma$,
$$ g \cdot \epsilon \left( \alpha\left(g\right)\right) = g = \epsilon \left(\beta \left(g\right)\right)\cdot g , \ \ \ g^{-1} \cdot g = \epsilon \left(\alpha\left(g\right)\right) , \ \ \ g \cdot g^{-1} = \epsilon \left(\beta\left(g\right)\right),$$
where we are denoting by $g^{-1}$ the image of $g$ by $i$. These maps will be called \textit{structure maps} and the groupoid will be denoted by $ \Gamma \rightrightarrows M$. We will also denote by $\Gamma_{x}$ (resp. $\Gamma^{x}$) the \textit{$\alpha$-fibre} $\alpha^{-1} \left( x \right)$ (resp. the \textit{$\beta$-fibre} $\beta^{-1} \left( x \right)$) and by $\Gamma_{x}^{y}$ the intersection $\alpha^{-1} \left( x \right) \cap \beta^{-1} \left( y \right)$. For each $ x \in M$, the group $\Gamma_{x}^{x}$ is called \textit{isotropy group at $x$}.\\
A groupoid $\Gamma \rightrightarrows M$ is said to be \textit{transitive} if the \textit{anchor map} $\left( \alpha , \beta \right): \Gamma \rightarrow M \times M$ is surjective. Equivalently, the sets $\Gamma_{x}^{y}$ are not empty for every $x,y \in M$. Notice that, in this case, all the isotropy groups are conjugated.\\
A \textit{groupoid morphism} between two groupoids, $\Gamma_{1} \rightrightarrows M_{1}$ and $\Gamma_{2} \rightrightarrows M_{2}$, consists of two maps $\Phi : \Gamma_{1} \rightarrow \Gamma_{2}$ and $\phi : M_{1} \rightarrow M_{2}$ such that for any $g_{1} \in \Gamma_{1}$
\begin{equation}\label{4}
\alpha_{2} \left( \Phi \left(g_{1}\right)\right) = \phi \left(\alpha_{1} \left(g_{1} \right)\right), \ \ \ \ \ \ \ \beta_{2} \left( \Phi \left(g_{1}\right)\right) = \phi \left(\beta_{1} \left(g_{1} \right)\right),
\end{equation}
where $\alpha_{i}$ and $\beta_{i}$ are the source and the target map of $\Gamma_{i} \rightrightarrows M_{i}$ respectively, for $i=1,2$, and preserves the composition, i.e.,
$$\Phi \left( g_{1} \cdot h_{1} \right) = \Phi \left(g_{1}\right) \cdot \Phi \left(h_{1}\right), \ \forall \left(g_{1} , h_{1} \right) \in \Gamma_{\left(2\right)}.$$
Observe that, as a consequence, $\Phi$ preserves the identities. We will denote this morphism as $\Phi$ (because, using equations \ref{4}, $\phi$ is completely determined by $\Phi$).\\
A \textit{Lie groupoid} is a groupoid $ \Gamma \rightrightarrows M$ where $\Gamma$ and $M$ are manifolds, the structure maps are differentiable and the source and the target maps are submersions. It is remarkable that, in the case of Lie groupoids, the sets $\Gamma_{x}$, $\Gamma^{x}$ and $\Gamma_{x}^{y}$ are manifolds (indeed, submanifolds of $\Gamma$). In fact, the isotropy groups are Lie groups.\\
A \textit{Lie groupoid morphism} is a morphism between Lie groupoids which is differentiable.\\
So, we define a \textit{(Lie) subgroupoid} of a (Lie) groupoid $\Gamma \rightrightarrows M$ as a (Lie) groupoid $\Gamma' \rightrightarrows M'$ such that $M' \subseteq M$, $\Gamma' \subseteq \Gamma$ and the inclusion maps induce a morphism of (Lie) groupoids. A \textit{reduction} of a transitive Lie groupoid is a transitive Lie subgroupoid over the same base.\\
\begin{remarkth}
\rm
There is a more abstract way of defining a groupoid. We can say that a groupoid is a ``small" category (the class of objects and the class of morphisms are sets) in which each morphism is invertible.

If $ \Gamma \rightrightarrows M$ is the groupoid, then $M$ is the set of objects and $\Gamma$ is the set of morphisms.

A groupoid morphism is a functor between these categories, which is a more natural definition.
\end{remarkth}

As in the case of Lie groups, we can define the concept of translation on a Lie groupoid: Let $\Gamma \rightrightarrows M$ be a groupoid. Then, for each $g \in \Gamma$ we may define the \textit{left translation on} $g$ as the map $L_{g} : \Gamma^{\alpha\left(g\right)} \rightarrow \Gamma^{\beta\left(g\right)}$ such that
$$ L_{g} \left( h \right) = g \cdot h,$$
for all $h \in  \Gamma^{\alpha\left(g\right)}$. The \textit{right translation} $R_{g}$ can be defined in a similar way. Obviously, both maps are diffeomorphisms with inverses $L_{g^{-1}}$ and $R_{g^{-1}}$ respectively.\\

\begin{example}[\textbf{Lie group}]
\rm
A Lie group $G$ is a Lie groupoid over a single point.
\end{example}

\begin{example}[\textbf{Pair groupoid}]
\rm
Let $M$ be a manifold. The product $M \times M$ is a Lie groupoid over $M$ such that the composition law is given by
$$\left( y,z \right) \cdot \left( x,y \right) = \left( x , z \right),$$
for all $\left( y,z \right), \left( x,y \right) \in M \times M$. This groupoid is called the \textit{pair groupoid on $M$}.\\
Notice that, for any (Lie) groupoid $\Gamma \rightrightarrows M$ the anchor map $\left( \alpha , \beta \right)$ is a morphism of (Lie) groupoids from $\Gamma$ to the pair groupoid on $M$.
\end{example}
\begin{example}[\textbf{Trivial groupoid}]
\rm
Let $M$ be a manifold and $G$ be a Lie group. As a natural generalization of the previous two examples (the Lie group and the pair groupoid) we have the groupoid $M \times M \times G \rightrightarrows M$. In this case, the composition is induced by the operation of the group. In fact, for each two elements $\left( y,z , g\right), \left( x,y , h \right) \in M \times M \times G$ we have that
$$\left( y,z, g \right) \cdot \left( x,y , h\right) = \left( x , z ,g\cdot h\right),$$
This groupoid is called the \textit{trivial groupoid on $M$ with group $G$}.
\end{example}

Next, we will present the \textit{$1-$jets groupoid} which will be one of the most important objects in what follows.\\
\begin{example}[\textbf{1-jets groupoid}]
\rm
Let $M$ be a manifold and $\Pi^{1} \left( M , M \right)$ be the set of all $1-$jets $j^{1}_{x,y} \phi$ of local diffeomorphisms $\phi: U \rightarrow V$ on $M$. Then, $\Pi^{1} \left( M , M \right)$ can be considered as a Lie groupoid over $M$ with the composition of $1-$jets as the composition law of the groupoid. Let $\left(x^{i}\right)$ and $\left( y^{j} \right)$ be local coordinates defined on two open subsets $U$ and $V$ of the base $M$ respectively, then we induce local coordinates on $\Pi^{1} \left( M , M \right)$ as follows
\begin{equation}\label{17}
\Pi^{1}\left(U,V\right) : \left(x^{i} , y^{j}, y^{j}_{i}\right),
\end{equation}
where, for each $ j^{1}_{x,y} \psi \in \Pi^{1}\left(U,V\right)$
\begin{itemize}
\item $x^{i} \left(j^{1}_{x,y} \psi\right) = x^{i} \left(x\right)$.
\item $y^{j} \left(j^{1}_{x,y}\psi \right) = y^{j} \left( y\right)$.
\item $y^{j}_{i}\left( j^{1}_{x,y}\psi\right)  = \left. \dfrac{\partial \left(y^{j}\circ \psi\right)}{\partial x^{i}}\right | _x$.
\end{itemize}
This groupoid will be called \textit{$1-$jets groupoid of $M$}.\\\\
\end{example}

\subsection*{Lie algebroids}
Now, we will introduce the notion of \textit{Lie algebroid}. Lie algebroids can be seen as a generalization of the Lie algebras and, as in the case of the Lie group, every Lie groupoid can be ``infinitesimally" described as a Lie algebroid.\\
A \textit{Lie algebroid} over a manifold $M$ is a triple $\left(  A \rightarrow M, \sharp , [ \cdot , \cdot ] \right)$, where $\pi : A \rightarrow M$ is a vector bundle together with a vector bundle morphism $\sharp : A \rightarrow TM$, called the \textit{anchor}, and a Lie bracket $[ \cdot  , \cdot ]$ on the space of sections, such that the Leibniz rule holds
\begin{equation}\label{18}
[ \alpha , f \beta ] = f [\alpha , \beta ] + \alpha^{\sharp}\left(f\right)\beta,
\end{equation}
for all $\alpha , \beta \in \Gamma \left(A\right)$ and $f \in \mathcal{C}^{\infty} \left(M\right)$. We are denoting ${\sharp}\left(\alpha\right)$ by $\alpha^{\sharp}$. $A$ is \textit{transitive} if $\sharp$ is surjective.\\ 
One can now prove that the anchor map is a morphism of Lie algebras, i.e.
\begin{equation}\label{20}
[\alpha , \beta]^{\sharp} = [\alpha^{\sharp}, \beta^{\sharp}] , \ \forall \alpha , \beta \in \Gamma \left(A\right).
\end{equation}
Notice that the bracket in the right side is just the Lie bracket of vector fields. We will use the same notation for both brackets.\\
Next, we should introduce the definition of a Lie algebroid morphism. However, the picture here is not so straightforward. The main problem is that a morphism between vector bundles does not, in general, induce a map between the modules of sections, so it is not immediately clear what should be meant by bracket relation.  We will give a direct definition in terms of $\left(\Phi,\phi\right)-$decompositons of sections which is easy to use, and is amenable to categorical methods.

Let $\left(A \rightarrow M , \sharp, [\cdot ,\cdot ]\right)$, $\left(A' \rightarrow M' , \sharp', [\cdot,\cdot]'\right)$ be Lie algebroids. Consider a vector bundle morphism $\Phi : A' \rightarrow A$, $\phi : M' \rightarrow M$ between $\pi : A \rightarrow M$ and $\pi' : A' \rightarrow M'$. We know that for each $ \alpha',\beta' \in \Gamma \left(A'\right)$, there exists $f_{i},g_{j} \in \mathcal{C}^{\infty} \left(M'\right)$ and $\alpha_{i}, \beta_{j} \in \Gamma \left(A\right)$ such that
$$\Phi \circ \alpha' = \sum_{i=1}^{k} f_{i} \left(\alpha_{i} \circ \phi \right), \ \ \ \sum_{j=1}^{k} g_{j}\left(\beta_{j} \circ \phi\right).$$
Thus, $\Phi$ is a morphism of Lie algebroids if it satisfies that
\begin{itemize}
\item[]\begin{equation}\label{27}
\sharp \circ \Phi = T\phi \circ \sharp',
\end{equation}

\item[]\begin{small}
\begin{equation}\label{28}
\Phi \circ [\alpha' , \beta' ] = \sum_{i,j=1}^{k}f_{i}g_{j}\left([ \alpha_{i} , \beta_{j} ] \circ \phi \right) + \sum_{j=1}^{k}{\alpha'}^{\sharp'}\left(g_{j}\right)\left(\beta_{j}\circ \phi\right) - \sum_{i=1}^{k}{\beta'}^{\sharp'}\left(f_{i}\right)\left(\alpha_{i}\circ \phi\right).
\end{equation}
\end{small}
\end{itemize}
In fact, the right-hand side of Eq. (\ref{28}) is independent of the choice of the $\left(\Phi,\phi\right)-$decompositions of $\alpha'$ and $\beta'$.\\
It is easy to prove that the composition preserves Lie agebroid morphisms and, hence, we can define the category of Lie algebroids.

\begin{remarkth}\label{29}
\rm
In particular, if $\alpha' \sim_{\left(\Phi , \phi\right)} \alpha$ and  $\beta' \sim_{\left(\Phi , \phi\right)} \beta$, then Eq. (\ref{28}) reduces to
$$\Phi \circ [\alpha', \beta'] = [\alpha , \beta]\circ \phi.$$

%Now consider two morphisms of Lie algebroids, $\Phi': A^{''} \rightarrow A'$, $\phi' : M^{''} \rightarrow M'$ and $\Phi: A' \rightarrow A$, $\phi : M' \rightarrow M$. One can observe that a $(\Phi', \phi'\right)-$decompositon,
%$$\Phi' \circ \alpha^{''} = \sum_{i=1}^{k} f_{i}^{''}(\alpha_{i}' \circ \phi'\right),$$
%combines with a $(\Phi , \phi\right)-$decompositon of each $\alpha_{i}'$ to give a $(\Phi \circ \Phi', \phi \circ \phi'\right)-$  decompositon, and verifies (\ref{28}) for decompositons so formed. Also, the condition on the anchors is easy to check. We therefore have a category of Lie algebroids. We will denote this category by $\mathcal{L}\mathcal{A}$.\\
%%\begin{remarkth}\label{29}
%%\rm
On the other hand, if $M=M'$ and $\phi = Id_{M}$ then Eq. (\ref{28}) reduces to 
$$ \Phi \circ [\alpha' , \beta'] = [\Phi \circ \alpha' , \Phi \circ \beta' ], \ \forall \alpha' , \beta' \in \Gamma \left(A'\right).$$ 
\end{remarkth}

A \textit{Lie subalgebroid} $A'$ of a Lie algebroid $A$ is a vector subbundle such that the inclusion is a  morphism of Lie algebroids. A Lie subalgebroid $A'$ of a transitive Lie algebroid $A$ is said to be a \textit{reduction of} $A$ if it is transitive and the base manifolds are equal.\\
Suppose that $M' \subseteq M$ is a closed submanifold then, using the $(i_{A'},i_{M'})-$ decomposition and extending functions, it satisfies that for all $\alpha' \in \Gamma (A')$ there exists $\alpha \in \Gamma (A)$ such that
$$ i_{A'} \circ \alpha' = \alpha \circ i_{M'}.$$ 
So, Eq. (\ref{28}) reduces to 
$$ i_{A'} \circ [\alpha' , \beta']_{M'} = [ \alpha ,  \beta ]_{M} \circ i_{M'}, \ \forall \alpha' , \beta' \in \Gamma (A').$$ 
As a third step, let us give a sketch of the construction of the Lie algebroid of a Lie groupoid. Let $\Gamma \rightrightarrows M$ be a Lie groupoid. Consider the vector bundle $\epsilon^{*}\left(Ker\left(T\beta\right)\right) : A\Gamma \rightarrow M$. Thus, the fibres will be interpreted as the tangent space of the $\beta-$fibres at the identities.\\
Now, consider the family $\frak X_{L} \left( \Gamma \right)$ of \textit{left-invariant vector fields on $\Gamma$}, i.e., the vector fields $\Theta$ on $\Gamma$ tangent to the $\beta-$fibres such that
\begin{equation}\label{100}
T_{h} L_{g} \left( \Theta \left( h \right) \right) = \Theta \left( g \cdot h \right),
\end{equation}
for all $g,h \in \Gamma$ with $\alpha \left( g \right) = \beta \left( h \right)$. Notice that, Eq. (\ref{100}) implies that the left-invariant vector fields are characterized by their image at the identities. With this, and using in fact that the left-invariant vector fields are tangent to the $\beta-$fibres, the space of sections of $A\Gamma$ is isomorphic to the space of left-invariant vector fields on $\Gamma$ and this fact gives us a Lie algebra structure over the space of section of $A \Gamma$. It is important to remark that the space of left-invariant vector fields on $\Gamma$ is closed under the Lie bracket of vector fields. Finally, we will construct the anchor $\sharp$ by restricting the tangent induced map of $\alpha$, i.e.,
$$\sharp \left( v_{x} \right) = T_{\epsilon \left( x \right) } \alpha \left( v_{x} \right),$$
for all $v_{x} \in A \Gamma_{x}$. It is not hard to prove that this three objects form a Lie algebroid called the \textit{Lie algebroid of the Lie groupoid} $\Gamma \rightrightarrows M$ and denoted by $A \Gamma$.\\
It is remarkable that each morphism between Lie groupoids induces a morphism between their respective Lie algebroids: mathematically speaking, this will imply that the construction induces a functor from the category of Lie groupoids to the category of Lie algebroids. In fact, let $\Phi : \Gamma \rightarrow \Gamma'$ be a Lie groupoid morphism from the groupoid $\Gamma \rightrightarrows M$ to $\Gamma' \rightrightarrows M$ over the identity map on $M$. Then, we can induce a map $A \Phi$ from $A \Gamma $ to $A \Gamma'$ in the following way:
\begin{equation}\label{101}
A \Phi \left(v_{x} \right) = T_{\epsilon \left( x \right)} \Phi_{x} \left( v_{x} \right),
\end{equation}
for all $v_{x} \in T_{\epsilon \left( x \right)} \Gamma^{x}$ where $ \Phi_{x}$ is the restriction of $\Phi$ to the $\beta-$fibre $\Gamma^{x}$. Then, $A \Phi$ is a Lie algebroid morphism from $A \Gamma$ to $A \Gamma'$.\\
In particular, a Lie subgroupoid $\overline{\Gamma} \rightrightarrows \overline{M}$ of a Lie groupoid $\Gamma \rightrightarrows M$ induces a Lie subalgebroid $A \overline{\Gamma}$ of a Lie algebroid $ A\Gamma$.\\
It is easy to realize that this construction generalizes the known one for Lie groups. So, if the Lie groupoid is a Lie group $G$, the resulting Lie algebroid is the associated Lie algebra $\frak g$. On the other hand, for a manifold $M$, it is immediate to check that the Lie algebroid of the pair groupoid $M \times M \rightrightarrows M$ is the tangent bundle $TM$ (with the structure of Lie algebroid induced by the Lie bracket of vector fields).\\
Finally, let $M \times M \times G$ be the trivial Lie groupoid on $M$ with group $G$. Then, the associated Lie algebroid can be interpreted as the vector bundle $A = TM \oplus \left( M \times \frak g\right) \rightarrow M$ such that
\begin{itemize}
\item[(i)] The anchor $\sharp :TM \oplus \left( M \times \frak g \right) \rightarrow TM$ is the projection.

\item[(ii)] Lie algebra structure over the space of sections is given by:
$$ \left[ \Theta \oplus f , \Xi \oplus g \right] =\left[ \Theta, \Xi \right] \oplus \{ \Theta \left( g \right) - \Xi \left( f \right)  + \left[f,g\right] \},$$
for all $  \Theta \oplus f , \Xi \oplus g \in \Gamma \left(A \right)$.
\end{itemize}
This Lie algebroid is called the \textit{trivial Lie algebroid on $M$ with structure algebra $\frak g$}.\\
Notice that, for each Lie algebroid $\left(A \rightarrow M , \sharp, [\cdot ,\cdot ]\right)$ the anchor $\sharp$ is a morphism of Lie algebroids from $A$ to the tangent algebroid $TM$.\\
Let us consider a Lie groupoid $\Gamma \rightrightarrows M$ and $A \Gamma$ its associated Lie algebroid. Then, the anchor $\sharp$ of $A \Gamma$ is the induced map of the anchor $\left( \alpha , \beta \right)$ of $\Gamma$, i.e.,
$$ A \left( \alpha , \beta \right)  =\sharp.$$
So, $\Gamma$ is transitive implies that $A \Gamma$ is transitive. In fact, the converse is also true and, hence, $\Gamma$ is transitive if, and only if, $A \Gamma$ is transitive.\\

Next, the \textit{$1-$jets algebroid} of a manifold $M$ will be the Lie algebroid $A \Pi^{1} \left( M , M \right)$ of the $1-$jets groupoid $\Pi^{1} \left( M , M \right)$ of $M$. We will give a detailed description below.\\
Let $\left(x^{i}\right)$ be a local coordinate system defined on some open subset $U \subseteq M$. Then, by using Eq. (\ref{101}), we can consider local coordinates on $A \Pi^{1}\left(M,M\right)$ as follows
\begin{equation}\label{42}
A\Pi^{1}\left(U,U\right) : \left(\left(x^{i} , x^{i}, \delta^{i}_{j}\right), 0 , v^{i}, v^{i}_{j}\right) \cong \left(x^{i} , v^{i}, v^{i}_{j}\right).
\end{equation}

Let $\Theta$ be a vector field on $M$. Denote by $\varphi^{\Theta}_{t}: U_{t} \rightarrow U_{-t}$ the (local) flow of $\Theta$. Then, for each $t$ we can construct a diffeomorphism,
$$ \Pi \varphi^{\Theta}_{t} : \Pi^{1} \left( U_{-t} , \mathcal{B} \right) \rightarrow \Pi^{1} \left( U_{t} , \mathcal{B} \right),$$
such that
$$ \Pi  \varphi^{\Theta}_{t}  \left( g \right) = g \cdot j^{1}_{\varphi^{\Theta}_{-t} \left( \alpha \left( g \right) \right),  \alpha \left( g \right) } \varphi^{\Theta}_{t}.$$
So, this flow induces a left-invariant vector field on $\Pi^{1} \left( M,M \right)$ which generates a section of $A \Pi^{1}\left( M,M \right)$ denoted by $j^{1}\Theta$. $j^{1}\Theta$ is called the \textit{complete lift of $\Theta$ on $\Pi^{1} \left( M,M \right)$}.\\
Let $\left( x^{i} \right)$ be a local chart of $M$ and $\left( x^{i} , y^{j} , y^{j}_{i} \right)$ be the induced local chart of $\Pi^{1} \left( M,M \right)$. Assume that, locally, $\Theta$ is written as follows,
$$\Theta = \Theta^{i}\dfrac{\partial}{\partial x^{i}}.$$
Then, locally, $j^{1}\Theta$ is expressed in the following way:
\begin{equation}\label{32}
j^{1}\Theta = - \Theta^{i}\dfrac{\partial}{\partial x^{i}} + \dfrac{\partial \Theta^{j}}{\partial x^{i}}\dfrac{\partial}{\partial y^{j}_{i}}.
\end{equation}
Notice that $j^{1}\Theta$ can be equivalently induced by a $1-$jet of $\Theta$. Thus, $A \Pi^{1} \left( M , M \right)$ can be interpreted as the bundle of $1-$jets of vector fields on $M$.\\
As a last step, we will introduce a new Lie algebroid to give a another interpretation of the $1-$jets Lie algebroid.\\

\begin{definition}
\rm
Let $ M$ be a manifold. A \textit{derivation on $M$} is a  $\mathbb{R}-$linear map $D : \frak X \left( M \right) \rightarrow \frak X \left( M \right)$ with a vector field $\Theta \in \frak X \left(M\right)$ such that for each $f \in \mathcal{C}^{\infty}\left(M\right)$ and $\Xi \in \frak X \left( M \right)$,
$$ D \left(f \Xi\right) = f D \left( \Xi \right) + \Theta\left(f\right) \Xi.$$
We call $\Theta$ the \textit{base vector field of $D$}. So, a derivation on $M$ is characterized by two geometrical objects, $D$ and $\Theta$.
\end{definition}
A classical example of derivation is given by the bracket of vector fields on a manifold $M$. In fact, let $\Theta$ be a vector field on $M$, the operator given by fixing $\Theta$ in the Lie bracket 
$$\left[\Theta , \cdot \right] : \frak X \left( M \right) \rightarrow \frak X \left( M \right),$$
is a derivation on $M$ with $\Theta$ as base vector field.\\
Another example comes from the so called \textit{covariant derivatives}. A covariant derivative on $M$ is a $\mathbb{R}-$bilinear map $\nabla : \frak X \left(M\right) \times \frak X \left(M \right) \rightarrow \frak X \left( M \right)$ such that,
\begin{itemize}
\item[(1)] It is $\mathcal{C}^{\infty} \left(M\right)-$linear in the first variable.
\item[(2)] For all $\Theta, \Xi \in \frak X \left( M \right)$ and $f \in \mathcal{C}^{\infty} \left(M\right)$,
\begin{equation}
\nabla_{\Theta} f \Xi = f \nabla_{\Theta}\Xi + \Theta \left( f\right) \Xi.
\end{equation}
\end{itemize}

Then, any vector field $\Theta \in \frak X \left(M\right)$ generates a derivation on $A$, $\nabla_{\Theta}$, (with base vector field $\Theta$) fixing the first coordinate again, i.e.,
$$ \nabla_{\Theta} : \frak X \left( M \right) \rightarrow \frak X \left( M \right),$$
such that
$$ \nabla_{\Theta}\left( \Xi \right) = \nabla_{\Theta} \Xi, \ \forall \Xi \in \frak X \left( M \right).$$
Associated to any covariant derivative $\nabla$ there are two important tensors:
\begin{itemize}
\item \textbf{Torsion:} $T \left( \Theta , \Xi \right)  = \nabla_{\Theta}\Xi - \nabla_{\Xi}\Theta - \left[\Theta,\Xi\right], \ \forall \Theta , \Xi \in \frak X \left( M \right).$

\item \textbf{Curvature:} $R \left( \Theta , \Xi \right)\chi =  \nabla_{\Theta} \nabla_{\Xi}\chi - \nabla_{\Xi} \nabla_{\Theta}\chi -\nabla_{\left[\Theta,\Xi\right]}\chi, \ \forall \Theta , \Xi ,\chi \in \frak X \left( M \right).$
\end{itemize}

A covariant derivative is said to be \textit{flat} if its curvature is zero.\\
We will use the well-known result:
\begin{lemma}\label{104}
Let $\nabla$ be a covariant derivative on a manifold $M$. $\nabla$ is flat and torsion-free if, and only if, there exists an atlas $\left( x^{i} \right)$ of $M$ such that
$$\nabla_{\mbox{\tiny $\dfrac{\partial}{\partial x^{j}}$}}\dfrac{\partial}{\partial x^{i}} = 0.$$

\end{lemma}
In general, for each local coordinates $\left( x^{i} \right)$ on $M$,
$$\nabla_{\mbox{\tiny $\dfrac{\partial}{\partial x^{j}}$}}\dfrac{\partial}{\partial x^{i}}   =\sum_{k} \Gamma^{k}_{i,j}  \dfrac{\partial}{\partial x^{k}}.$$
The local functions $\Gamma^{k}_{i,j}$ are called \textit{Christoffel symbols of $\nabla$ respect to} $\left( x^{i} \right)$.\\

Now, the space of derivations on $M$ can be considered as the space of sections of a vector bundle $\mathfrak{D}\left(TM\right)$ on $M$. We can endow this vector bundle with a Lie algebroid structure.
\begin{itemize}
\item Let $D_{1} , D_{2}$ be derivations on $M$, we can define $[D_{1} , D_{2}]$ as the commutator, i.e.,
$$[D_{1} , D_{2}] = D_{1} \circ D_{2} - D_{2}\circ D_{1}.$$
A simple computation shows that the commutator of two derivations is again a derivation, indeed, the base vector field of $[D_{1} , D_{2}]$ is given by 
\begin{equation}
[\Theta_{1} , \Theta_{2}],
\end{equation}
where $\Theta_{1}$ and $\Theta_{2}$ are the base vector fields of $D_{1}$ and $D_{2}$ respectively.
\item Let $D$ be a derivation on $M$, then $D^{\sharp}$ is its base vector field.
\end{itemize}
Thus, with this structure $\mathfrak{D}\left(TM\right)$ is a transitive Lie algebroid called the \textit{Lie algebroid of derivations on $M$}. The space of sections of $\mathfrak{D}\left(TM\right)$, the derivations on $M$, will be denoted by $Der\left(TM\right)$.\\
Note that in this Lie algebroid the fibre-wise linear sections of $\sharp$ are $\mathcal{C}^{\infty}\left(M\right)-$linear maps from $\frak X \left(M\right)$ to $Der\left(TM\right)$. So, the space of fibre-wise linear sections of $\sharp$ is, indeed, the space of covariant derivatives on $M$. In fact, it is easy to see that a covariant derivative $\nabla$ is a Lie algebroid morphism (from the tangent algebroid to the algebroid of derivations) if, and only if, $\nabla$ is flat.\\
Finally, it is turn to relate this algebroid with the $1-$jets Lie algebroid. Consider $\Lambda \in \Gamma \left( A\Pi^{1} \left( M , M \right)\right)$ and $\Theta^{\Lambda}$ its associated left-invariant vector field on $\Pi^{1} \left( M , M \right)$. Denote by $\Phi_{t}^{\Lambda}: \mathcal{U}_{t} \rightarrow \mathcal{U}_{-t}$ the flow of $\Theta^{\Lambda}$.\\
Then, we can define a (local) linear map $\left(\Phi_{t}^{\Lambda}\right)^{*} : \frak X \left( M \right) \rightarrow \frak X \left( M \right)$ satisfying
$$ \{\left(\Phi_{t}^{\Lambda}\right)^{*} \left( \Theta\right) \} \left( x \right)  =  \Phi_{t}^{\Lambda} \left(\epsilon\left( x\right)\right) \left( \Theta \left(\left( \alpha \circ \Phi_{t}^{\Lambda}\right)\left(\epsilon \left(x\right)\right)\right)\right),$$
for each $\Theta \in \frak X \left( M \right)$ and $x \in M$. 
Thus, we can define the following derivation on $M$,
$$ D^{\Lambda}  =\left. \dfrac{\partial}{\partial t}  \left(\Phi_{t}^{\Lambda}\right)^{*}\right | _{t=0}.$$
In other words, for each $ \Theta \in \frak X \left( M \right)$ and $x \in M$ we have
\begin{eqnarray*}
D^{\Lambda}\Theta \left( x \right) &=& \left. \dfrac{\partial}{\partial t} \left( \Phi_{t}^{\Lambda} \left(\epsilon\left( x\right)\right) \left( \Theta \left(\left( \alpha \circ \Phi_{t}^{\Lambda}\right)\left(\epsilon \left(x\right)\right)\right)\right)\right)  \right | _{t=0}.
\end{eqnarray*}
Notice that, for all $f \in \mathcal{C}^{\infty} \left( M \right)$
\begin{small}
\begin{eqnarray*}
D^{\Lambda}f \Theta \left( x \right) &=& \left. \dfrac{\partial}{\partial t} \left( \Phi_{t}^{\Lambda}  \left( \epsilon \left( x \right) \right)   \left( f\left( \left( \alpha \circ \Phi_{t}^{\Lambda} \right) \left(\epsilon \left(x \right) \right) \right)\Theta \left(\left( \alpha \circ \Phi_{t}^{\Lambda} \right) \left( \epsilon \left( x \right) \right)\right) \right)\right)\right | _{t=0}\\
&=& \left. \dfrac{\partial}{\partial t} \left( f\left(\left( \alpha \circ \Phi_{t}^{\Lambda} \right) \left(\epsilon \left( x \right) \right)\right) \Phi_{t}^{\Lambda}  \left(\epsilon \left( x \right) \right)   \left( \Theta \left(\left( \alpha \circ \Phi_{t}^{\Lambda} \right) \left(\epsilon \left( x \right) \right)\right)\right)\right)\right | _{t=0}\\
 &=& \Lambda^{\sharp} \left( x \right) \left( f \right) \Theta \left( x \right) + f \left( x \right) D^{\Lambda}\Theta \left( x \right).
\end{eqnarray*}

\end{small}
It is immediate to prove that for each $\Theta \in \frak X \left( M \right)$ one has that
\begin{equation}\label{102}
D^{j^{1}\Theta}\Xi = \left[\Theta,\Xi\right], \ \forall \Xi \in \frak X \left( M \right).
\end{equation}

This construction gives us a linear map between the sections of the $1-$jets Lie algebroid and the algebroid of derivations which induces a Lie algebroid isomorphism $\mathcal{D} : A \Pi^{1} \left( M , M \right) \rightarrow  \mathfrak{D} \left(A\right)$ over the identity map on $M$.\\
Notice that using this isomorphism, we can consider a one-to-one map from fibre-wise linear sections of $\sharp$ in $ A \Pi^{1}\left(M,M\right)$ to covariant derivatives over $M$. Thus, having a fibre-wise linear section $\Delta$ of $\sharp$ in \linebreak $A \Pi^{1}\left(M,M\right)$ we will denote its associated covariant derivative by $\nabla^{\Delta}$. Furthermore, $\Delta$ is a morphism of Lie algebroids if, and only if, $\nabla^{\Delta}$ is flat.\\
Next, we will show how locally the map $\mathcal{D}$ looks more natural:

\begin{lemma}\label{76}
Let $M$ be a manifold and $\Lambda $ be a section of the $1-$jets algebroid with local expression
$$\Lambda \left(x^{i}\right) = \left(x^{i} , \Lambda^{j} , \Lambda^{j}_{i}\right).$$
The matrix $\Lambda^{j}_{i}$ is (locally) the associated matrix to $D^{\Lambda}$, i.e.,
$$D^{\Lambda} \left( \dfrac{\partial}{\partial x^{i}} \right)= \sum_{j} \Lambda^{j}_{i}  \dfrac{\partial}{\partial x^{j}},$$
and the base vector field of $D^{\Lambda}$ is $\Lambda^{\sharp}$ which is given locally by $\left(x^{i} ,  \Lambda^{j} \right)$.
\end{lemma}

Let $\Delta$ be a fibre-wise linear section of $\sharp$ in $A \Pi^{1}\left(M,M\right)$ and $\nabla^{\Delta}$ be its associated covariant derivative. Thus, for each $\left(x^{i}\right)$ local coordinate system on $M$
$$\Delta\left(x^{i} , \dfrac{\partial}{\partial x^{j}}\right) = \left(x^{i} , \dfrac{\partial}{\partial x^{j}}, \Delta^{j}_{i}\right),$$
where $\Delta^{j}_{i}$ depends on $\dfrac{\partial}{\partial x^{j}}$ linearly. Thus, we will change the notation as follows
\begin{equation}\label{145}
\Delta^{j}_{i} \left(x^{l} , \dfrac{\partial}{\partial x^{k}}\right) = \Delta^{j}_{i,k} \left( x^{l} \right).
\end{equation}
Therefore,
$$ \nabla^{\Delta}_{\dfrac{\partial}{\partial x^{j}}} \dfrac{\partial}{\partial x^{i}} = D^{\Delta \left( \dfrac{\partial}{\partial x^{j}}\right)} \dfrac{\partial}{\partial x^{i}} =\sum_{k} \Delta^{k}_{i,j}  \dfrac{\partial}{\partial x^{k}},$$
where $\Delta \left( \dfrac{\partial}{\partial x^{j}}\right)$ is the (local) section of $A\Pi^{1}\left(M,M\right)$ given by
$$ \Delta \left( \dfrac{\partial}{\partial x^{j}}\right) \left(x\right) = \Delta\left(x\right) \left( \left. \dfrac{\partial}{\partial x^{j}} \right | _{x}\right).$$

So, $\Delta^{k}_{i,j}$ are just the Christoffel symbols of $\nabla^{\Delta}$.\\

%Using this remark, an section of $\sharp$ in the $1-$jets Lie algebroid, $\Lambda$,  is integrable if, and only if, there exists an atlas on $M$, $\mathcal{H}$, such that $\nabla^{\Lambda}$ is a covariant derivative with Christoffel symbols equal to zero, i.e.,
%$$\nabla^{\Lambda}_{ \dfrac{\partial}{\partial x^{i}} } \dfrac{\partial}{\partial x^{j}} = 0, \ \forall i,j.$$

\section{Uniformity and Homogeneity}

A \textit{body} $\mathcal{B}$ is a three-dimensional differentiable manifold which can be covered with just one chart. An embedding $\phi : \mathcal{B} \rightarrow \mathbb{R}^{3}$ is called a \textit{configuration of} $\mathcal{B}$ and its $1-$jet $j_{X,\phi \left(X\right)}^{1} \phi$ at $X \in \mathcal{B}$ is called an \textit{infinitesimal configuration at $X$}. We usually identify the body with any one of its configurations, say $\phi_{0}$, called \textit{reference configuration}. Given any arbitrary configuration $\phi$, the change of configurations $\kappa = \phi \circ \phi_{0}^{-1}$ is called a \textit{deformation}, and its $1-$jet $j_{\phi_{0}\left(X\right) , \phi \left(X\right)}^{1} \kappa$ is called an \textit{infinitesimal deformation at $\phi_{0}\left(X\right)$}.

From now on we make the following identification: $\mathcal{B} \cong \phi_{0} \left( \mathcal{B}\right)$.

For elastic bodies, the mechanical response of a material is completely characterized by one function $W$ which depends, at each point $X \in \mathcal{B}$, on the gradient of the deformations evaluated at the point. Thus, $W$ is defined (see \cite{MELZA}) as a differentiable map denoted by the same letter
$$ W : Gl\left(3 , \mathbb{R}\right) \times \mathcal{B} \rightarrow V,$$
where $V$ is a real vector space. Another equivalent way of considering $W$ is as a differentiable map
$$ W : \Pi^{1} \left( \mathcal{B}, \mathcal{B}\right) \rightarrow V,$$
which does not depend on the final point, i.e., for all $X,Y,Z \in \mathcal{B}$
\begin{equation}\label{92}
 W \left( j_{X,Y}^{1} \phi\right) = W \left( j_{X,Z}^{1} \left( \tau_{Z-Y} \circ \phi\right)\right), \ \forall j_{X,Y}^{1}\phi \in \Pi^{1} \left( \mathcal{B}, \mathcal{B}\right),
\end{equation}
where $\tau_{v}$ is the translation map on $\mathbb{R}^{3}$ by the vector $v$. This map will be called the \textit{response functional}.

Notice that, using Eq. (\ref{92}), we can define $W$ over $\Pi^{1} \left( \mathcal{B} , \mathbb{R}^{3}\right)$, which is the open subset of $\Pi^{1} \left( \mathbb{R}^{3} , \mathbb{R}^{3}\right)$ given by the $1-$jets of local diffeomorphisms from points of $\mathcal{B}$ to points of $\mathbb{R}^{3}$.\\

Now, suppose that an infinitesimal neighbourhood of the material around the point $Y$ can be grafted so perfecly into a neighbourhood of $X$, that the graft cannot be detected by any mechanical experiment. If this condition is satisfied with every point $X$ of $ \mathcal{B}$, the body is said \textit{uniform}. We can express this physical property in a geometric way as follows.

\begin{definition}
\rm
A body $\mathcal{B}$ is said to be \textit{uniform} if for each two points $X,Y \in \mathcal{B}$ there exists a local diffeomorphism $\psi$ from an open neighbourhood $U \subseteq \mathcal{B}$ of $X$ to an open neighbourhood $V \subseteq \mathcal{B}$ of $Y$ such that $\psi \left(X\right) =Y$ and
\begin{equation}\label{50}
W \left( j^{1}_{Y, \kappa \left(Y\right)} \kappa \cdot j^{1}_{X,Y} \psi \right) = W \left( j^{1}_{Y, \kappa \left(Y\right)} \kappa\right),
\end{equation}
for all infinitesimal deformation $j^{1}_{Y , \kappa \left(Y\right)} \kappa$. The $1-$jets of local diffeomorphisms satisfying Eq. (\ref{50}) are called \textit{material isomorphisms}.
\end{definition}
These kind of maps are going to be important. Let us show that the set of these maps can be endowed of a groupoid structure. For each two points we will denote by $G \left(X,Y\right)$ the collection of all $1-$jets $j_{X,Y}^{1}\psi$ which satisfy Eq. (\ref{50}). So, the set of all material ismorphisms can be written as follows,
$$\Omega \left( \mathcal{B}\right) = \cup_{X,Y \in \mathcal{B}} G\left(X,Y\right).$$
Notice that the identities are material isomorphisms, and the composition and the inversion of $1-$jets preserve Eq. (\ref{50}). Hence, $\Omega \left( \mathcal{B} \right)$ has structure of groupoid over $\mathcal{B}$ which is, indeed, a subgroupoid of the $1-$jets groupoid $\Pi^{1}\left( \mathcal{B} , \mathcal{B}\right)$. We will denote $\alpha^{-1} \left(X\right)$ (resp. $\beta^{-1} \left(X\right)$) by $\Omega_{X}\left( \mathcal{B}\right)$ (resp. $\Omega^{X}\left( \mathcal{B}\right)$). The elements of the isotropy group $G \left( X , X \right)$ will be called \textit{material symmetries at $X$}.

So, the following result is obvious.
\begin{proposition}
Let $\mathcal{B}$ be a body. $\mathcal{B}$ is uniform if, and only if, $\Omega \left( \mathcal{B}\right)$ is a transitive subgroupoid of $\Pi^{1} \left( \mathcal{B} , \mathcal{B}\right)$.
\end{proposition}
Notice that, at general, we cannot ensure that $\Omega \left( \mathcal{B}\right) \subseteq \Pi^{1} \left( \mathcal{B} , \mathcal{B} \right)$ is a Lie subgroupoid (see for instance \cite{MGEOEPS,MD,CHARDIST}). Our assumption is that $\Omega \left( \mathcal{B}\right)$ is in fact a Lie subgroupoid and, in this case, $\Omega \left(\mathcal{B}\right)$ is said to be the \textit{material groupoid of $\mathcal{B}$}.\\
As we have seen, a body is uniform is the function $W$ depends on the point $X$ precisely according to Equation (\ref{50}). In addition, a body is said to be \textit{homogeneous} if we can choose a global section of the material groupoid which is constant on the body, more precisely:

\begin{definition}\label{51}
\rm
A body $\mathcal{B}$ is said to be \textit{homogeneous} if it admits a global deformation $\kappa$ which induces a global section of $\left(\alpha , \beta\right)$ in $\Omega \left( \mathcal{B}\right)$, $P$, i.e., for each $X,Y \in \mathcal{B}$
$$ P\left(X, Y\right) = j^{1}_{X,Y} \left(\kappa^{-1} \circ \tau_{\kappa\left(Y\right) - \kappa \left(X\right)} \circ \kappa\right),$$
where $\tau_{\kappa\left(Y\right) - \kappa \left(X\right)}: \mathbb{R}^{3} \rightarrow \mathbb{R}^{3}$ denotes the translation on $\mathbb{R}^{3}$ by the vector $\kappa\left(Y\right) - \kappa \left(X\right)$. $\mathcal{B}$ is said to be \textit{locally homogeneous} if there exists a covering of $\mathcal{B}$ by homogeneous open sets.
\end{definition}
Suppose that $\mathcal{B}$ is homogeneous. Then, if we take global coordinates $\left(x^{i}\right)$ given by the induced diffeomorphism $\kappa$, we deduce that $P$ is locally expressed by
\begin{equation}\label{52}
P\left(x^{i},y^{j}\right) = \left(x^{i},y^{j} , \delta^{j}_{i}\right).
\end{equation}
If $\mathcal{B}$ is locally homogeneous we can cover $\mathcal{B}$ by local sections of $\left(\alpha , \beta\right)$ in $\Omega \left( \mathcal{B}\right)$ which satisfy Eq. (\ref{52}). The (local) coordinates generated by these $\kappa's$ will be called \textit{homogeneous coordinates}.

\section{Integrability}
In this section we will do a purely mathematic development to study the notion of integrability of reduced subgroupoids of the 1-jets\linebreak groupoid. This notion of integrability will be closely related with the notion of integrability of $G-$structures.

Note that there exists a Lie groupoid isomorphism $L : \Pi^{1}\left( \mathbb{R}^{n},\mathbb{R}^{n}\right) \rightarrow \mathbb{R}^{n}\times \mathbb{R}^{n} \times Gl\left(n , \mathbb{R} \right)$ over the identity map defined by
$$ L \left( j^{1}_{x,y}\phi \right) = \left(x,y, j^{1}_{0,0}\left(\tau_{-y} \circ \phi \circ \tau_{x}\right)\right), \ \forall j^{1}_{x,y}\phi \in \Pi^{1}\left(\mathbb{R}^{n},\mathbb{R}^{n}\right),$$
where $\tau_{v}$ denote the translation on $\mathbb{R}^{n}$ by the vector $v \in \mathbb{R}^{n}$ and we are identifying $Gl \left( n , \mathbb{R}\right)$ with the space of linear isomorphisms from $T_{0} \mathbb{R}^{n}$ to $T_{0} \mathbb{R}^{n}$ (i.e. the fibre of the frame bundle $F \mathbb{R}^{n}$ at $0$). So, if $G$ is a Lie subgroup of $Gl\left(n, \mathbb{R}\right)$, we can transport $\mathbb{R}^{n} \times \mathbb{R}^{n} \times G$ by this isomorphism to obtain a reduced Lie subgroupoid of $\Pi^{1}\left( \mathbb{R}^{n} , \mathbb{R}^{n}\right)$. This kind of reduced subgroupoids will be called \textit{standard flat} on $\Pi^{1}\left(\mathbb{R}^{n},\mathbb{R}^{n}\right)$.

Let $U, V \subseteq M$ be two open subsets of $M$. We denote by $\Pi^{1}\left(U,V\right)$ the open subset of $\Pi^{1}\left(M,M\right)$ defined by $\left(\alpha , \beta \right)^{-1}\left(U \times V\right)$. Note that if $U=V$, then, $\Pi^{1}\left(U, U \right)$ is, indeed, the 1-jets groupoid of $U$ and, in this way, our notation is consistent.\\

\begin{definition}
\rm
A reduced subgroupoid $\Pi^{1}_{G}\left(M,M\right)$ of $\Pi^{1}\left(M,M\right)$ will be called \textit{integrable} if it is locally diffeomorphic to a standard flat on $\Pi^{1}\left(\mathbb{R}^{n},\mathbb{R}^{n}\right)$
\end{definition}

Before continuing with our discussion, we need to explain what we understand by ``locally diffeomorphic" in this case. $\Pi^{1}_{G}\left(M,M\right)$ is locally diffeomorphic to $\mathbb{R}^{n} \times \mathbb{R}^{n} \times G \rightrightarrows \mathbb{R}^{n}$ for some Lie subgroup $G$ of $Gl\left(n, \mathbb{R}\right)$ if for all $x, y\in M$ there exist two open sets $U, V \subseteq M$ with $x \in U$, $y \in V$ and two local charts, $\psi_{U}: U \rightarrow \overline{U}$ and $\psi_{V}: V \rightarrow \overline{V}$, which induce a diffeomorphism
\begin{equation}\label{53}
 \Psi_{U,V} : \Pi^{1}_{G}\left(U,V\right) \rightarrow \overline{U} \times \overline{V} \times G,
\end{equation}
such that $ \Psi_{U, V} = \left( \psi_{U} \circ \alpha, \psi_{V} \circ \beta , \overline{\Psi}_{U,V}\right)$, where
$$ \overline{\Psi}_{U,V} \left( j_{x,y}^{1} \phi \right) = j_{0,0}^{1} \left( \tau_{-\psi_{V} \left(y\right)} \circ \psi_{V} \circ \phi \circ \psi_{U}^{-1} \circ \tau_{\psi_{U}\left(x\right)}\right), \ \forall  j_{x,y}^{1} \phi \in \Pi^{1} \left(U,V\right).$$
Notice that, $\Pi^{1}_{G}\left(U,V\right)$ and $ \overline{U} \times \overline{V} \times G$ are Lie groupoids if, and only if, $U = V$ and $ \overline{U} = \overline{V}$. Suppose that $U=V$ and $\overline{U} = \overline{V}$, then, for all $x \in U$ $\Psi_{U,U}\left(j_{x,x}^{1}Id\right) \in G$. However, $\Psi_{U,U}\left(j_{x,x}^{1}Id\right)$ is not necessarily the identity map and, hence, $\Psi_{U,U}$ is not an isomorphism of Lie groupoids.\\

\begin{proposition}\label{54}
Let $\Pi^{1}_{G} \left(M,M\right)$ be a reduced Lie subgroupoid of \linebreak$\Pi^{1}\left(M,M\right)$. $\Pi_{G}\left(M,M\right)$ is integrable if, and only if, we can cover $M$ by local charts $\left(\psi_{U} , U\right)$ which induce Lie groupoid isomorphisms from $\Pi^{1}_{G} \left(U,U\right)$ to the restrictions of the standard flat over $G$ to $\psi_{U}\left(U\right)$.\\

\begin{proof}

On the one hand, suppose that $\Pi^{1}\left(M,M\right)$ is integrable. Let $x_{0} \in M$ be a point in $M$ and $\psi_{U}: U \rightarrow \overline{U}$ and $\psi_{V}: V \rightarrow \overline{V}$ be local charts through $x_{0} $ which induced diffeomorphism
$$\Psi_{U, V} : \Pi^{1}_{G}\left(U,V\right)  \rightarrow \overline{U} \times \overline{V} \times G.$$
For each $y \in U \cap V$,
$$ \overline{\Psi}_{U,V} \left( j_{y,y}^{1} Id \right) = j_{0,0}^{1} \left( \tau_{-\psi_{V} \left(y\right)} \circ \psi_{V} \circ \psi_{U}^{-1} \circ \tau_{\psi_{U}\left(y\right)}\right) \in G.$$
Then, for all $j_{x,y}^{1} \phi \in \Pi^{1}_{G} \left(U \cap V , U \cap V\right)$, we have
\begin{small}
$$ j_{0,0}^{1} \left( \tau_{-\psi_{U} \left(y\right)} \circ \psi_{U} \circ \phi \circ \psi_{U}^{-1} \circ \tau_{\psi_{U}\left(x\right)}\right) =   $$
$$ j_{0,0}^{1} \left( \tau_{-\psi_{U} \left(y\right)} \circ \psi_{U} \circ \psi_{V}^{-1} \circ \tau_{\psi_{V}\left(y\right)}\right) \cdot j_{0,0}^{1} \left( \tau_{-\psi_{V} \left(y\right)} \circ \psi_{V} \circ \phi \circ \psi_{U}^{-1} \circ \tau_{\psi_{U}\left(x\right)}\right) \in G.$$
\end{small}
Therefore, denoting $U \cap V$ by $W$, the map
$$\Psi_{W,W} : \Pi^{1}_{G}\left(W,W\right) \rightarrow \overline{W} \times \overline{W} \times G,$$
is, indeed, a Lie groupoid isomorphism over $\psi_{W}$ where \linebreak $ \Psi_{W,W} = \left( \psi_{W} \circ \alpha, \psi_{W} \circ \beta , \overline{\Psi}_{W,W}\right)$, $\psi_{W}$ is the restriction of $\psi_{U}$ to $W$ and for all $j_{x,y}^{1} \phi \in \Pi^{1}\left(W,W\right)$,
$$ \overline{\Psi}_{W,W} \left( j_{x,y}^{1} \phi \right) = j_{0,0}^{1} \left( \tau_{-\psi_{W} \left(y\right)} \circ \psi_{W} \circ \phi \circ \psi_{W}^{-1} \circ \tau_{\psi_{W}\left(y\right)}\right).$$
On the other hand, suppose that for each $x \in M$ there exists a local chart $\left( \psi_{U} , U\right)$ through $x$ which induces a Lie groupoid isomorphism over $\psi_{U}$, namely
\begin{equation}\label{55}
\Psi_{U, U} : \Pi^{1}_{G}\left(U,U\right) \rightarrow \overline{U} \times \overline{U} \times G,
\end{equation}
such that $ \Psi_{U,U} = \left( \psi_{U} \circ \alpha, \psi_{U} \circ \beta , \overline{\Psi}_{U,U}\right)$, where for each \linebreak $ j_{x,y}^{1} \phi \in \Pi^{1} \left(U,U\right)$,
$$ \overline{\Psi}_{U,U} \left( j_{x,y}^{1} \phi \right) = j_{0,0}^{1} \left( \tau_{-\psi_{U} \left(y\right)} \circ \psi_{U} \circ \phi \circ \psi_{U}^{-1} \circ \tau_{\psi_{U}\left(x\right)}\right).$$
Take open sets $U,V \subseteq M$ such that there exist $\psi_{U}$ and $\psi_{V}$ satisfy Eq. (\ref{55}). 
Suppose that $ U \cap V \neq \emptyset$. Then, for all $x,y \in U \cap V$, we have
\begin{equation}\label{56}
 j_{0,0}^{1} \left( \tau_{-\psi_{U}\left(y\right)} \circ \psi_{U} \circ \psi_{V}^{-1} \circ \tau_{\psi_{V}\left(y\right)}\right) \cdot j_{0,0}^{1} \left( \tau_{-\psi_{V}\left(x\right)} \circ \psi_{V} \circ \psi_{U}^{-1} \circ \tau_{\psi_{U}\left(x\right)}\right) \in G.
\end{equation}

Fixing $Z \in U \cap V$, we consider
$$ j_{0,0}^{1} \left( \tau_{-\psi_{U}\left(Z\right)} \circ \psi_{U} \circ \psi_{V}^{-1} \circ \tau_{\psi_{V}\left(Z\right)}\right) \triangleq  A \in Gl\left(n, \mathbb{R}\right).$$
Furthermore, it is also true that
\begin{equation}\label{conjugatedsubgroup3224}
 A \cdot G \cdot A^{-1} = G.
\end{equation}
In fact,
$$ G = \{ j_{0,0}^{1} \left( \tau_{-\psi_{U} \left(y\right)} \circ \psi_{U} \circ \phi \circ \psi_{U}^{-1} \circ \tau_{\psi_{U}\left(y\right)}\right) \ / \ j_{y,y}^{1}\phi \in \Pi^{1}_{G}\left(U,U\right)\}.$$
So, we define the diffeomorphism $\overline{\psi}_{V} \triangleq A \cdot \psi_{V} : V \rightarrow A \cdot \overline{V}$. Then, using Eq. (\ref{56}) for all $ y \in U \cap V$, we deduce that
$$ j_{0,0}^{1} \left( \tau_{-\overline{\psi}_{V}\left(y\right)} \circ \overline{\psi}_{V} \circ \psi_{U}^{-1} \circ \tau_{\psi_{U}\left(y\right)}\right) =  j_{0,0}^{1} A\cdot \left( \tau_{-\psi_{V}\left(y\right)} \circ \psi_{V} \circ \psi_{U}^{-1} \circ \tau_{\psi_{U}\left(y\right)}\right) $$
\begin{equation}\label{57}
 =  A \cdot j_{0,0}^{1} \left( \tau_{-\psi_{V}\left(y\right)} \circ \psi_{V} \circ \psi_{U}^{-1} \circ \tau_{\psi_{U}\left(y\right)}\right) \in G.
\end{equation}
In this way, we consider
$$
\begin{array}{rccl}
\Psi_{U,V} : & \Pi^{1}_{G}\left(U,V\right) & \rightarrow & \overline{U} \times A \cdot \overline{V} \times G \\
& j^{1}_{x,y}\phi  &\mapsto &  \left(\psi_{U}\left(x\right), \overline{\psi}_{V}  \left(y\right), \overline{\Psi}_{U,V}\left(j^{1}_{x,y} \phi \right)\right).
\end{array}
$$
where,
$$ \overline{\Psi}_{U,V}\left(j^{1}_{x,y} \phi \right) = j_{0,0}^{1} \left( \tau_{-\overline{\psi}_{V} \left(y\right)} \circ \overline{\psi}_{V} \circ \phi \circ \psi_{U}^{-1} \circ \tau_{\psi_{U}\left(x\right)}\right).$$
We will check that $\overline{\Psi}_{U,V}$ is well-defined. We fix $j_{x,y}^{1} \phi \in \Pi^{1}_{G} \left(U,V\right)$. Then, we can consider two cases:
\begin{itemize}
\item[(i)] $y \in U \cap V$. Then, using Eq. (\ref{57})
\begin{small}
$$ j_{0,0}^{1} \left( \tau_{-\overline{\psi}_{V} \left(y\right)} \circ \overline{\psi}_{V} \circ \phi \circ \psi_{U}^{-1} \circ \tau_{\psi_{U}\left(x\right)}\right) = $$
$$  j_{0,0}^{1} \left( \tau_{-\overline{\psi}_{V} \left(y\right)} \circ \overline{\psi}_{V} \circ  \psi_{U}^{-1} \circ \tau_{\psi_{U}\left(y\right)}\right) \cdot j_{0,0}^{1} \left( \tau_{-\psi_{U} \left(y\right)} \circ \psi_{U} \circ \phi \circ \psi_{U}^{-1} \circ \tau_{\psi_{U}\left(x\right)}\right) \in G.$$
\end{small}

\item[(ii)] $ y \notin U \cap V$. Then, 
$$ j_{z,x}^{1} \left( \psi^{-1}_{V} \circ \tau_{\psi_{V}\left(z\right) - \psi_{V} \left(y\right)} \circ \psi_{V} \circ \phi \right) \triangleq j_{z,x}^{1} \phi_{z} \in \Pi^{1}_{G}\left(M,M\right).$$
Hence,
\begin{small}
$$ j_{0,0}^{1} \left( \tau_{-\overline{\psi}_{V} \left(y\right)} \circ \overline{\psi}_{V} \circ \phi \circ \psi_{U}^{-1} \circ \tau_{\psi_{U}\left(x\right)}\right) =$$
$$ A \cdot j_{0,0}^{1} \left( \tau_{-\psi_{V} \left(y\right)} \circ \psi_{V} \circ \phi \circ \psi_{U}^{-1} \circ \tau_{\psi_{U}\left(x\right)}\right) = $$
$$A \cdot j_{0,0}^{1} \left( \tau_{-\psi_{V} \left(z\right)} \circ \psi_{V}  \circ \phi_{z} \circ \psi_{U}^{-1} \circ \tau_{\psi_{U}\left(x\right)}\right) = $$
$$j_{0,0}^{1} \left( \tau_{-\overline{\psi}_{V} \left(z\right)} \circ \overline{\psi}_{V}  \circ \phi_{z} \circ \psi_{U}^{-1} \circ \tau_{\psi_{U}\left(x\right)}\right) \in G.$$
\end{small}

\end{itemize}

Thus, it is immediate to prove that $\Psi_{U,V}$ is a diffeomorphism which commutes with the restrictions of the structure maps.

Finally, if $U \cap V = \emptyset$ we can find a finite family of local neighbourhoods $\{ V_{i} \}_{i=1,...,k}$ such that
\begin{itemize}
\item[(i)] $U = V_{1}$
\item[(ii)] $V= V_{k}$
\item[(iii)] $V_{i} \cap V_{i+1} \neq \emptyset, \ \forall i$
\end{itemize}
Thus, we can find $\Psi_{U,V}$ following a similar procedure than above.\\
\end{proof}
\end{proposition}
Using this result, we can prove that a reduced subgroupoid $\Pi_{G}^{1}\left(M,M\right)$ of $\Pi^{1}\left(M,M\right)$ is integrable respect to two different Lie groups $G$ and $\tilde{G}$ if, and only if, $G$ and $\tilde{G}$ are conjugated.\\
There is a special reduced subgroupoid of $\Pi^{1}\left(M,M\right)$ which will play an important role in the following. A trivial reduced subgroupoid of $\Pi^{1}\left(M,M\right)$ or \textit{parallelism of} $\Pi^{1}\left(M,M\right)$ is a reduced subgroupoid of $\Pi^{1}\left(M,M\right)$, $\Pi_{e}^{1}\left(M,M\right) \rightrightarrows M$, such that for each $x,y \in M$ there exists a unique $1-$jet $j^{1}_{x,y}\phi \in \Pi^{1}_{e}\left(M,M\right)$. So, equivalently, a parallelism can be interpret as a section $P: M \times M \rightarrow \Pi^{1}\left(M,M\right)$ of $\left(\alpha, \beta\right)$ which is a morphism of Lie groupoids (over the identity map on $M$). Using this, we can also speak about \textit{integrable sections of} $\left(\alpha , \beta\right)$.\\

Now, let us consider the induced coordinates given in Eq. (\ref{17}). Then, an integrable section $P$ can be written locally as follows,
$$P\left(x^{i}, y^{j} \right) = \left(x^{i} , y^{j} , \delta^{j}_{i}\right),$$
or equivalently
\begin{equation}\label{60}
P \left(x,y\right) = j^{1}_{x,y} \left(\psi^{-1} \circ \tau_{\psi \left(y\right) - \varphi \left(x\right)} \circ \varphi \right),
\end{equation}
for some two local charts $\left( \varphi , U\right) , \left(\psi , V\right)$ on $M$.

%%Notice that, using Proposition \ref{54}, $P$ is an integrable section if and only if we can cover $M$ by local charts $\left(\varphi , U\right)$ such that
%%\begin{equation}\label{61}
%%P_{|U} \left(x,y\right) = j^{1}_{x,y} \left(\varphi^{-1} \circ \tau_{\varphi \left(y\right) - \varphi \left(x\right)} \circ \varphi \right).
%%\end{equation}

Next, analogously to the case of $G-$structures, we can characterize the integrable subgroupoids using (local) integrable sections. However, in this case it is not so easy because having a reduced subgroupoid we do not have a fixed structured group $G$. Let $\Pi^{1}_{G}\left(M,M\right)$ be a reduced subgroupoid of $\Pi^{1}\left(M,M\right)$ and $Z_{0} \in FM$, where $FM$ is the linear frame bundle of $M$ (see \cite{DIGFMCORD} for a detailed study on the geometry of frame bundles), be a linear frame at $z_{0} \in M$. Then, we define
\begin{equation}\label{62}
G  := \{ Z_{0}^{-1} \cdot g\cdot Z_{0} \ / \ g \in \Pi^{1}_{G}\left(M,M\right)_{z_{0}}^{z_{0}}\} = Z_{0}^{-1}\cdot \Pi^{1}_{G}\left(M,M\right)_{z_{0}}^{z_{0}}\cdot Z_{0},
\end{equation}
Therefore, $G$ is a Lie subgroup of $Gl\left(n , \mathbb{R}\right)$. This Lie group will be called \textit{associated Lie group} to $\Pi^{1}_{G}\left(M,M\right)$. Note that, as a difference with $G-$structures, we don't have a unique Lie group $G$ but all of them are conjugated.
\begin{proposition}\label{63}
A reduced subgroupoid $\Pi^{1}_{G}\left(M,M\right)$ of $\Pi^{1}\left(M,M\right)$ is integrable if, and only if, for each two points $x ,y \in M$ there exist coordinate systems $\left(x^{i}\right)$ and $\left(y^{j}\right)$ over $U,V \subseteq M$, respectively with $x \in U$ and $y \in V$ such that the local section, 
\begin{equation}\label{64}
P\left(x^{i}, y^{j}\right) = \left(x^{i}, y^{j} , \delta^{j}_{i}\right),
\end{equation}
takes values into $\Pi^{1}_{G}\left(M,M\right)$.
\begin{proof}
First, it is obvious that if $\Pi^{1}_{G}\left(M,M\right)$ is integrable then, we can restrict the maps $\Psi_{U,V}^{-1}$ to $\overline{U} \times \overline{V} \times \{e\}$ to get (local) integrable sections of $\left(\alpha, \beta\right)$ which takes values on $\Pi^{1}_{G}\left(M,M\right)$.

Conversely, in a similar way to Proposition \ref{54} we claim that for each $x \in M$ there exists an open set $U \subseteq M$ with $x \in U$ and $P : U \times U \rightarrow  \Pi^{1}_{G} \left(U,U\right)$ an integrable sections of $\left(\alpha , \beta\right)$ given by
$$ P\left(x,y\right) = j_{x,y}^{1}\left(\psi_{U}^{-1} \circ \tau_{\psi_{U}\left(y\right) - \psi_{U}\left(x\right)} \circ \psi_{U}\right),$$
where $\psi_{U} : U \rightarrow \overline{U}$ is a local chart at $x$.

Then, we can build the map
$$\Psi_{U,U}^{-1} : \overline{U} \times \overline{U} \times \{e\} \rightarrow \Pi^{1}_{G}\left(U,U\right),$$
defined in the obvious way.

Now, let $z_{0} \in U$ be a point at $U$, $Z_{0} \triangleq j_{0,z_{0}}^{1}\left(\psi_{U}^{-1} \circ \tau_{\psi_{U}\left(z_{0}\right)}\right) \in FU$ be a frame at $z_{0}$ and $G$ be the Lie subgroup satisfying Eq. (\ref{62}). Then, we can define
$$\Psi_{U, U} : \Pi^{1}_{G}\left(U,U\right) \rightarrow \overline{U} \times \overline{U} \times G,$$
where for each $j_{z_{0},z_{0}}^{1}\phi \in \Pi^{1}_{G} \left(z_{0}\right)$ and $x,y \in \overline{U}$ we define
$$ \Psi_{U,U}^{-1}\left( x,y , Z^{-1}_{0}\cdot j_{z_{0},z_{0}}^{1}\phi \cdot Z_{0}\right) $$
$$ = j_{0,\psi_{U}^{-1}\left(y\right)}^{1}\left(\psi_{U}^{-1} \circ \tau_{y}\right) \cdot  [ Z^{-1}_{0} \cdot j_{z_{0},z_{0}}^{1}\phi \cdot Z_{0} ] \cdot  j_{\psi_{U}^{-1}\left(x\right),0}^{1} \left( \tau_{- x} \circ \psi_{U}\right).$$
Hence the map $\Psi_{U,U}: \overline{U} \times \overline{U} \times G \rightarrow \Pi^{1}_{G}\left(U,U\right)$ is an isomorphism of Lie groupoids induced by $\psi_{U}$.

To end the proof, we only have to use Proposition \ref{54}.
\end{proof}
\end{proposition}

%%%%
As we have noticed, this notion of integrability is closely related with the integrability of $G-$structures. In fact, for each $G-$structure $\omega_{G} \left( M \right)$ on $M$ we consider the following set,
\begin{equation}\label{103}
\mathcal{G}\left(\omega_{G}\left(M\right)\right) = \{ L_{y} \cdot [L_{x}^{-1}]\ / \ L_{x}, L_{y} \in \omega_{G}\left(M\right)\}.
\end{equation}
It is not hard to prove that $\mathcal{G}\left(\omega_{G}\left(M\right)\right)$ is a reduced subgroupoid of $\Pi^{1}\left(M, M\right)$. Conversely, for all reduced subgroupoid $\Pi^{1}_{G}\left(M, M\right)$ of $\Pi^{1}\left(M, M\right)$ and all frame $Z_{0} \in FM$ at $z_{0} \in M$, there exists a $G-$structure,
$$\omega_{G}\left(M\right)= \Pi^{1}_{G}\left(M, M\right)_{z_{0}} \cdot Z_{0},$$
such that
$$ \mathcal{G}\left(\omega_{G}\left(M\right)\right) = \Pi^{1}_{G}\left(M, M\right).$$
With this, we can claim that a $G-$structure $\omega_{G}\left(M\right)$ is integrable if, and only if, its associated Lie subgroupoid $ \mathcal{G}\left(\omega_{G}\left(M\right)\right)$ of $\Pi^{1} \left( M , M \right)$ is integrable.

\vspace{1cm}

Now, we want to work with the notion of integrability in the associated Lie algebroid of the 1-jets groupoid. Let $U \subseteq M$ be an open subset of $M$. We denote by $A\Pi^{1}\left(U,U\right)$ the open Lie subalgebroid of $A\Pi^{1}\left(M,M\right)$ defined by the associated Lie algebroid of $\Pi^{1}\left(U,U\right)$.\\

\begin{definition}\label{66}
\rm
Let $A\Pi^{1}_{G}\left(M,M\right)$ be a Lie subalgebroid of $A \Pi^{1} \left(M,M\right)$. $A\Pi^{1}_{G}\left(M,M\right)$ is said to be \textit{integrable by $G$} if it is locally isomorphic to the trivial algebroid $T\mathbb{R}^{n} \oplus \left( \mathbb{R}^{n} \times  \mathfrak{g}\right)$, where $\mathfrak{g}$ is the Lie algebra of the Lie subgroup $G$ of $Gl\left(n ,\mathbb{R}\right)$.
\end{definition}

Again, ``locally isomorphic" means that we are inducing local coordinates from the base in the natural way. It is not hard to prove that a Lie subalgebroid $A\Pi^{1}_{G}\left(M,M\right)$ of $A \Pi^{1} \left(M,M\right)$ is integrable by $G$ if, and only if, it is the associated Lie algebroid of an integrable Lie subgroupoid $\Pi^{1}_{G}\left(M,M\right)$ of $\Pi^{1}\left(M,M\right)$.\\

Analogously to the case of $1-$jets groupoid, a \textit{parallelism} of \linebreak$A\Pi^{1}\left(M,M\right)$ is an associated Lie algebroid of a parallelism of $\Pi^{1}\left(M,M\right)$. Hence, a parallelism is a section of $\sharp$ which is a morphism of Lie algeborids. In this way, we will also speak about \textit{integrable sections of $\sharp$}. Notice that, any morphism of Lie algebroid $\Delta : TM \rightarrow A \Pi^{1} \left( M , M \right)$ can be integrated to a morphism of Lie groupoids from $ P: M \times M \rightarrow \Pi^{1} \left( M , M\right)$ such that
$$ A P = \Delta.$$
This is just an application of a generalization of the Lie's second fundamental theorem from Lie groups to Lie groupoids (see for instance \cite{IMJMO}).\\

Let $\left(x^{i}\right)$ be a local coordinate system defined on some open subset $U \subseteq M$. Notice that each integrable section of $\left(\alpha , \beta\right)$ in $\Pi^{1}\left(M,M\right)$, $P$, is a Lie groupoid morphism. Hence, $P$ induces a Lie algebroid morphism $AP : TM \rightarrow A\Pi^{1}\left(M,M\right)$ which is a section of $\sharp$. So, taking into account that, locally, 
$$ P \left(x^{i} , y^{j} \right) = \left(x^{i} , y^{j} , \delta^{j}_{i}\right),$$
we have that each integrable section can be written locally as follows
$$ A P \left(x^{i}, \dfrac{\partial}{\partial x^{i}} \right) = \left(x^{i} , \dfrac{\partial}{\partial x^{i}} , 0\right).$$

Now, using Proposition \ref{63}, we have the following analogous proposition.
\begin{proposition}\label{70}
A reduced subalgebroid $A\Pi^{1}_{G}\left(M,M\right)$ of $A\Pi^{1}\left(M,M\right)$ is integrable by $G$ if, and only if, there exist local integrable sections of $\sharp$ covering $M$ which takes values on $A\Pi^{1}_{G}\left(M,M\right)$.
\end{proposition}
Equivalently, for each point $ x \in M$ there exists a local coordinate system $\left(x^{i}\right)$ over an open set $U \subseteq M$ with $x \in U$ such that the local sections
$$\Delta \left(x^{i}, \dfrac{\partial}{\partial x^{i}} \right) = \left(x^{i} , \dfrac{\partial}{\partial x^{i}} , 0\right),$$
are in $A \Pi^{1}_{G}\left(M,M\right)$.\\

Finally, we will use the algebroid of derivations on $T M$. Thus, taking into account that the map $\mathcal{D} : \Gamma \left(A \Pi^{1}\left( M,M\right)\right) \rightarrow Der\left(T M\right)$ defines a Lie algebroid isomorphism $\mathcal{D} : A \Pi^{1}\left( M,M\right) \rightarrow  \mathfrak{D} \left(T M\right)$ over the identity map on $M$ and using Lemma \ref{104} we can give another characterization of the integrability over the $1-$jets algebroid.

\begin{proposition}

Let $\Delta$ be a fibre-wise linear section of $\sharp$ in the $1-$jets Lie algebroid, $A \Pi^{1} \left(M,M\right)$. Then, it is integrable if, and only if, the covariant derivative $\nabla^{\Delta}$is flat and torsion-free. Equivalently, a $\Delta$ section of $\sharp$ which is a morphism of Lie algebroids is integrable if, and only if, the covariant derivative $\nabla^{\Delta}$ is torsion-free.

\end{proposition}

This permits us to say that a reduced subgroupoid $\Pi^{1}_{G}\left(M,M\right)$ of $\Pi^{1}\left(M,M\right)$ is integrable if, and only if, $M$ can be covered by (local) torsion-free flat covariant derivatives which takes vales in \linebreak$\mathcal{D} \left( A\Pi^{1}_{G}\left(M,M\right) \right)$.\\

\section{Application to material bodies}
This section is devoted to apply our results to the theory of material bodies.\\
Let $\mathcal{B}$ be a body. Taking into account the definition of homogeneity (see Definition \ref{51}) and the above result we can state the following proposition:

\begin{proposition}\label{65}
Let $\mathcal{B}$ be a uniform body. If $\mathcal{B}$ is homogeneous then $\Omega \left( \mathcal{B} \right)$ is integrable. Conversely, $\Omega \left( \mathcal{B} \right)$ is integrable implies that $\mathcal{B}$ is locally homogeneous.
\end{proposition}

Notice that Eq. (\ref{103}) proves that our definition of homogeneity (which is given using the material groupoid) is, indeed, equivalent to that used in \cite{MELZA} (see \cite{MELZASEG} or \cite{CCWANSEG}; see also \cite{FBLOO} and \cite{GAMAU2}) where the authors use $G-$structures to characterize this property.\\

Next, let us consider the induced subalgebroid of material groupoid, $A \Omega \left( \mathcal{B} \right)$. This Lie algebroid will be called \textit{ material algebroid of $\mathcal{B}$}.\\

Take $\Lambda \in \Gamma  \left( A\Omega \left( \mathcal{B}\right)\right)$. So, the flow of the left-invariant vector field $X_{\Lambda}$, $\{ \varphi^{\Lambda}_{t} \}$, can be restricted to $\Omega\left( \mathcal{B} \right)$.\\
Hence, for any infinitesimal deformation $ g $, we have
$$ W\left( \varphi^{\Lambda}_{t} \left( g \right) \right) = W\left( \varphi^{\Lambda}_{t} \left( g \cdot \epsilon \left( \alpha \left(g \right) \right) \right) \right) = W\left( g \cdot \varphi^{\Lambda}_{t} \left(  \epsilon \left( \alpha \left(g \right) \right) \right) \right) = W \left(g \right).$$
Thus, for each $g \in \Pi^{1} \left(  \mathcal{B} , \mathcal{B} \right)$, we deduce
$$T W \left( X_{\Lambda} \left( g \right)\right) = \left. \dfrac{\partial}{\partial t} \left(W \left(\varphi^{\Lambda}_{t}\left( g \right) \right)\right) \right | _{0}= \left. \dfrac{\partial}{\partial t} \left(W \left( g \right) \right) \right | _{0} = 0.$$
Therefore,
\begin{equation}\label{93}
T W \left( X_{\Lambda}\right) = 0.
\end{equation}

Conversely, it is easy to prove that Eq. (\ref{93}) implies that $ \Lambda \in \Gamma \left( A \Omega \left( \mathcal{B} \right)\right)$.

In this way, the material algebroid can be defined without using the material groupoid by imposing Eq. (\ref{93}). Thus, we can characterize the homogeneity and uniformity using only the material Lie algebroid.

\begin{proposition}\label{73}
Let $\mathcal{B}$ be a simple body. Then, $\mathcal{B}$ is uniform if, and only if, $A\Omega \left( \mathcal{B} \right)$ is transitive.

\end{proposition}
Let us assume that a section $\Lambda$ of the material algebroid $A \Omega \left( \mathcal{B} \right)$ has  the following local expression,
$$ \Lambda \left( x^{i}  \right) = \left(  x^{i}, \Lambda^{i} , \Lambda^{j}_{i}  \right),$$
for a given local coordinates $\left( x^{i}\right)$ on $\mathcal{B}$. Then, Eq. (\ref{93}) can be written as follows,
\begin{equation}\label{79}
\Lambda^{k} \dfrac{\partial  W}{\partial x^{k}} + F^{r}_{m} \Lambda^{m}_{l}\dfrac{\partial  W}{\partial x^{r}_{l}} = 0.
\end{equation}
So, the uniformity of the material is linked to the properties of Eq. (\ref{79}).
\begin{proposition}\label{74}
Let $\mathcal{B}$ be a uniform body. If $\mathcal{B}$ is homogeneous, then, $A\Omega \left( \mathcal{B} \right)$ is integrable by a Lie subgroup $G$ of $Gl\left(n, \mathbb{R}\right)$. Conversely, if $A\Omega \left( \mathcal{B} \right)$ is integrable by $G$ then $\mathcal{B}$ is locally homogeneous.

\end{proposition}

Using Proposition \ref{70}, this result can be expressed locally as follows.

\begin{proposition}\label{75}
Let $\mathcal{B}$ be a uniform body. $\mathcal{B}$ is locally homogeneous if and only if for each point $x  \in \mathcal{B}$ there exists a local coordinate system $\left(x^{i}\right)$ over $U \subseteq \mathcal{B}$ with $x \in U$ such that the local section of $\sharp$, 
$$\Delta \left(x^{i}, \dfrac{\partial}{\partial x^{i}} \right) = \left(x^{i} , \dfrac{\partial}{\partial x^{i}} , 0\right),$$
takes values in $A \Omega\left( \mathcal{B}\right)$.
\end{proposition}

Therefore, denoting by $ \mathcal{D}\left( \mathcal{B}\right)$ to the Lie subalgebroid of the derivation algebroid on $\mathcal{B}$, $ \mathcal{D} \left( A \Omega \left( \mathcal{B} \right) \right) \leq \mathfrak{D} \left(T \mathcal{B}\right)$, we can give the following result:

\begin{theorem}\label{78}
Let $\mathcal{B}$ be a uniform body.
\begin{itemize}
\item[(i)] $\mathcal{B}$ is homogeneous, if and only if, there exists a (global) torsion-free and flat covariant derivative on $\mathcal{B}$ which takes values in $\mathcal{D}\left( \mathcal{B}\right)$.

\item[(ii)] $\mathcal{B}$ is locally homogeneous if, and only if, $\mathcal{B}$ can be covered by (local) torsion-free and flat covariant derivatives which takes vales in $\mathcal{D} \left( \mathcal{B} \right)$.
\end{itemize} 
\end{theorem}

\section{Example}
We will use Proposition \ref{75} to work with an example. We will consider a model of a so-called \textit{simple liquid crystal.} These simple materials were introduced by Coleman \cite{COLE} and Wang \cite{CCWANTHIRD}.\\
%%%
%%%Notice that, taking into account Eq. (\ref{93}) and Proposition \ref{75}, a simple body $\mathcal{B}$ is locally homogeneous if, and only if, for all $x  \in \mathcal{B}$ there exists a local coordinate system $\left(x^{i}\right)$ on $x \in U$
%%%\begin{equation}\label{105}
%%%\dfrac{\partial}{\partial x^{i}} \left( W \right) = 0.
%%%\end{equation}

Let $\mathcal{B}$ be a simple body (we will assume that $\mathcal{B}$ is an open subset of $\mathbb{R}^{3}$ by taking the image by the reference configuration) with a mechanical response $ W : \Pi^{1} \left( \mathcal{B} , \mathcal{B} \right) \rightarrow V$ such that for all $ h = j^{1}_{X,Y} \phi \in \Pi^{1} \left( \mathcal{B} , \mathcal{B} \right)$ we have

$$ W \left( h \right) = \widehat{W} \left( r \left( h \right) , J \left( h \right)\right),$$
\noindent
where, denoting by $F$ the associated matrix to $j^{1}_{X,Y} \phi$ (with respect to the canonical basis of $\mathcal{B}$),
\begin{itemize}
\item $r \left( j^{1}_{X,Y} \phi \right) = g \left( Y \right) \left( T_{X} \phi \left(  e \left( X  \right)\right),  T_{X} \phi \left( e \left( X  \right)\right) \right)$
\item $ J \left( j^{1}_{X,Y} \phi \right) = det \left( F \right)$
\end{itemize}
with $e \in \frak X \left( \mathcal{B} \right)$ a vector field which is not zero at any point and $g$ a Riemannian metric on $\mathbb{R}^{3}$. Notice that the tangent bundle $T \mathcal{B}$ is canonically isomorphic to $ \mathcal{B} \times \mathbb{R}^{3}$. So, for each $Y \in \mathcal{B}$, $g \left( Y \right)$ can be seen as a inner product on $\mathbb{R}^{3}$. Then, the expression of $r$ turns into the following,
$$r \left( j^{1}_{X,Y} \phi \right) = g \left( Y \right) \left( F \cdot \left(  e^{i} \left( X  \right)\right),  F \cdot \left( e^{i} \left( X  \right)\right) \right),$$
where, by using the canonical isomorphism $T \mathcal{B} \rightarrow \mathcal{B} \times \mathbb{R}^{3}$, for all $X \in \mathcal{B}$ $e \left( X \right) = \left( X , e^{i} \left( X \right) \right)$. We will use both expressions with the same notation.\\
%%%Notice that we are trivializing the tangent bundle to write $g \left( X \right) \left( F\cdot  e \left( X  \right), F \cdot  e \left( X  \right) \right)$. Maybe, in a rigorous way, we could write $ T_{X} L_{F} \left( e \left( X \right)\right)$ where $L_{F} : \mathbb{R}^{3} \rightarrow \mathbb{R}^{3}$ is the diffeomorphism given by
%%%$$ L_{F} \left( v \right) = F \cdot v, \ \forall v \in \mathbb{R}^{3}.$$
%%%Let $\left( x^{i} \right)$ be a local system of coordinates in $\mathbb{R}^{3}$. Then, for all $X \in Dom \left( x^{i} \right)$
%%%$$ T_{X} L_{F} \left( \dfrac{\partial}{\partial x^{i}_{|X}} \right) = F \cdot \dfrac{\partial}{\partial x^{i}_{|F \cdot X}}.$$
%%%So, the tangent map of $L_{F}$ at a point $X$ applied to a vector $V$ is essentially $F \cdot V$.\\
Now, we want to study the condition which characterizes the material algebroid: A left-invariant vector field $\Theta \in \frak X_{L} \left( \Pi^{1} \left( \mathcal{B} , \mathcal{B} \right)\right)$ restricts to a section of $A \Omega \left( \mathcal{B} \right)$ if, and only if,
$$  \Theta \left( W\right) = 0.$$
So, we should study $TW$ over left-invariant vector fields. Let $\Theta \in \frak X_{L} \left( \Pi^{1} \left( \mathcal{B} , \mathcal{B} \right)\right)$ be a left-invariant vector field and consider the canonical local system of coordinates $\left( X^{i} \right)$ in $\mathbb{R}^{3}$. We will denote by $\left( X^{i} , Y^{j} , F^{j}_{i} \right)$ the induced local coordinates of $\left( X^{i} \right)$ in $\Pi^{1} \left( \mathcal{B} , \mathcal{B} \right)$. The local expression of $\Theta$ will be denoted as follows,
$$ \Theta \left( X^{i} , Y^{j} , F^{j}_{i} \right) = \left( \left( X^{i} , Y^{j} , F^{j}_{i} \right)  , \delta X^{i} ,0,  F^{j}_{l} \delta P^{l}_{i} \right).$$ 
Now, we will begin given the derivatives of $r$ and $J$. For each $W \in gl \left( 3, \mathbb{R} \right)$ and $v \in \mathbb{R}^{3}$ we have that,

\begin{eqnarray*}
&\textbf{(i)}&\ \ \ \left. \dfrac{\partial r}{\partial X} \left( v \right) \right | _{j_{X,Y}^{1}\phi} =  2   g \left( Y\right)\left( T_{X}\phi \left( e \left( X \right)\right) ,T_{X}\phi \left(  \left. \dfrac{\partial e}{\partial X}\left( v \right) \right | _{X}  \right)\right).\\
&\textbf{(ii)}&\ \ \ \left. \dfrac{\partial r}{\partial F} \left(  W  \right) \right | _{j_{X,Y}^{1}\phi}= 2   g \left( Y\right)\left(F \cdot \left(  e^{i} \left( X  \right)\right) ,W \cdot \left(e^{i} \left( X \right) \right) \right).\\
&\textbf{(iii)}&\ \ \ \left. \dfrac{\partial J}{\partial F} \left(  W  \right) \right | _{j_{X,Y}^{1}\phi}= det \left( F \right) Tr \left( F^{-1} \cdot W \right).
\end{eqnarray*}
Here $F$ is the Jacobian matrix of $\phi$ at $X$ and $ \left. \dfrac{\partial e }{\partial X}  \left( v \right)\right | _{X} $ is the vector at $X$ such that
$$  \left. \dfrac{\partial e }{\partial X} \left( v \right)\right | _{X} = \left( X , \left. \dfrac{\partial e^{i}}{\partial X^{l}} \right | _{X} v^{l} \right).$$
Hence, $\Theta$ restricts to a section of the material algebroid $A \Omega \left( \mathcal{B} \right)$ if, and only if,

\begin{eqnarray*}
 0 &=& 2 \left. \dfrac{\partial \widehat{W}}{\partial r} \right | _{j_{X,Y}^{1}\phi}  g \left( Y\right)\left( T_{X}\phi \left( e \left( X \right)\right) ,T_{X}\phi \left(   \left. \dfrac{\partial e}{\partial X} \left(\delta X^{i} \left(X\right)\right)\right | _{X}\right)\right)+\\
&&+ 2 \left. \dfrac{\partial \widehat{W}}{\partial r}\right | _{j_{X,Y}^{1}\phi}   g \left( Y\right)\left( F \cdot \left(  e^{i} \left( X  \right)\right) ,F^{j}_{l} \delta P^{l}_{i}\left(X \right) \cdot \left(e^{i} \left( X \right) \right) \right)+\\
&& + det \left(F \right) \left. \dfrac{\partial \widehat{W}}{\partial J} \right | _{j_{X,Y}^{1}\phi} Tr \left( \delta P^{j}_{i}\left(X\right)\right),
\end{eqnarray*}
for all $j_{X,Y}^{1} \phi \in \Pi^{1} \left( \mathcal{B} , \mathcal{B} \right)$. So, a sufficient but not necessary condition would be,

\begin{eqnarray*}
&\textbf{(1)}&\ \ \ Tr \left( \delta P^{j}_{i}\left(X\right)\right) = 0.\\
&\textbf{(2)} &\ \ \ \mbox{\scriptsize $g \left( Y\right)\left( F \cdot \left(  e^{i} \left( X  \right)\right)  , F^{r}_{m} \cdot \left(   \left. \dfrac{\partial e}{\partial X} \left(\delta X^{i}\left(X\right) \right) \right | _{X}+ \delta P^{j}_{i}  \left(X\right)\cdot \left(e^{i} \left( X \right)\right)  \right)\right) = 0$},
\end{eqnarray*}
for all $1-$jets of local diffeomorphisms $j_{X,Y}^{1} \phi $ on $ \mathcal{B}$. By using that $g$ is non-degenerate and $e \left( X \right)$ is non-zero, we turn these conditions into the following
\begin{eqnarray*}\label{72}
&\textbf{(1)'}&\ \ \   \delta P^{i}_{i} = 0.\\
&\textbf{(2)'}&\ \  \dfrac{\partial e^{j}}{\partial X^{l}}   \delta X^{l} + \delta P^{j}_{l}  e^{l} =0, \ \forall j,
\end{eqnarray*}
where $e^{j}$ are the coordinates of $e$ respect to $\left( X^{j} \right)$.\\
Let us now study the uniformity of the material. By using Proposition \ref{73} $\mathcal{B}$ is uniform if, and only if, the material algebroid of $\mathcal{B}$ is transitive.\\
Let $V_{X} = \left( X , V^{i}\right)$ be a vector at $X \in \mathcal{B}$. Then, we should find a (local) left-invariant vector field $ \Theta$ such that
\begin{itemize}
\item $\Theta \left( W\right) = 0.$

\item $ T_{\epsilon \left( X \right) } \alpha \left( \Theta \left( \epsilon \left( X \right) \right) \right) = V_{X},$
\end{itemize}
where $\epsilon$ and $\alpha$ are the identities map and the source map of the material groupoid respectively.\\
Let us fix the local expression of $\Theta$ as follows,
$$\Theta \left( X^{i} , Y^{j} , F^{j}_{i} \right) = \left( \left( X^{i} , Y^{j} , F^{j}_{i} \right)  , \delta X^{i} ,0,  F^{j}_{l} \delta P^{l}_{i} \right).$$
Then,
$$ T_{\epsilon \left( X \right) } \alpha \left( \Theta \left( \epsilon \left( X \right) \right) \right) = \left(  X^{i}  \left(X\right), \delta X^{i} \left( X \right) \right).$$
So, it should satisfy that,
$$  \delta X^{i} \left( X \right)  = V^{i}, \ \forall i.$$
By taking into account identities \textbf{(1)'} and \textbf{(2)'}, it is enough to find a family of (local) maps $A^{j}_{i}$ from the body to the space of matrices satisfying that
\begin{eqnarray*}\label{71}
&\textbf{(1)''}&\ \ \   A^{i}_{i}= 0.\\
&\textbf{(2)''}&\ \  \dfrac{\partial e^{j}}{\partial X^{l}}   V^{l} = -A^{j}_{l}  e^{l} , \ \forall j,
\end{eqnarray*}
It is just an easy exercise to prove that there are infinite solutions $A^{j}_{i}$ of the equations \textbf{(2)'} and \textbf{(2)'} and, hence, $\mathcal{B}$ is uniform.\\

From now on, we will assume that $\widehat{W}$ is an immersion. In that way, $\textbf{(1)'}$ and $\textbf{(2)'}$ are also necessary conditions.\\

Next, we will study the condition of (local) homogeneity. As we know (Proposition \ref{75}) $\mathcal{B}$ is (locally) homogeneous if, and only if, there exists a local system of coordinates $\left( x^{i} \right)$ such that the local section of $\sharp$,
$$\Delta \left( x^{i} , \dfrac{\partial }{\partial x^{i}} \right) = \left( x^{i} , \dfrac{\partial }{\partial x^{i}}, 0 \right),$$
takes values in the material algebroid $A \Omega \left( \mathcal{B} \right)$. Equivalently,
\begin{equation}\label{105}
\dfrac{\partial W}{\partial x^{i}} = 0 , \ \forall i.
\end{equation}
So, let us study this equality. Notice that,
$$\dfrac{\partial W}{\partial x^{i}} = \dfrac{\partial \widehat{W}}{\partial r} \dfrac{\partial r}{\partial x^{i}} + \dfrac{\partial \widehat{W}}{\partial J} \dfrac{\partial J}{\partial x^{i}}.$$

Thus, by using that $\widehat{W}$ is an immersion, $\left( x^{i} \right)$ are homogeneous coordinates if, and only if,

\begin{eqnarray*}
&\textbf{(1)'''}&\ \ \ \dfrac{\partial r}{\partial x^{i}} = 0, \ \forall i.\\
&\textbf{(2)'''}&\ \ \ \dfrac{\partial J}{\partial x^{i}} = 0, \ \forall i.
\end{eqnarray*}

Observe that the form of $\widehat{W}$ is not important to evaluate the (local) homogeneity of $\mathcal{B}$ as long as $\widehat{W}$ is an immersion.\\
Let $\left( x^{i} \right)$ be a system of homogeneous coordinates on $\mathcal{B}$. Then, for each $j_{X,Y}^{1} \phi \in \Pi^{1} \left( \mathcal{B} , \mathcal{B} \right)$
\begin{eqnarray*}
r \left( j_{X,Y}^{1} \phi \right) & = & g \left( Y\right)\left( T_{X}\phi \left( e \left( X \right)\right) ,  T_{X}\phi \left( e \left( X \right)\right) \right)\\
& = &  g \left( Y\right)\left( T_{X}\phi \left( e^{i} \left( X \right)\left. \dfrac{\partial }{\partial x^{i}}\right | _{X}\right) ,  T_{X}\phi \left( e^{j} \left( X \right)\left. \dfrac{\partial }{\partial x^{j}}\right | _{X}\right) \right)\\
& = & e^{i} \left( X \right)e^{j} \left( X \right) \left. \dfrac{\partial \phi^{k}}{\partial x^{i}}\right | _{X} \left. \dfrac{\partial \phi^{l} }{\partial x^{j}} \right | _{X} g_{kl} \left( Y \right),
\end{eqnarray*}
where, in this case, $e^{j}$ are the coordinates of $e$ respect to $\left(x^{i} \right)$. So, considering the induced coordinates $\left( x^{i} , y^{j} , y^{j}_{i}\right)$ of $\left( x^{i} \right)$ on $\Pi^{1} \left( \mathcal{B} , \mathcal{B} \right)$ we have that
\begin{eqnarray*}
r \circ \left( x^{i} , y^{j} , y^{j}_{i}\right)^{-1} \left( \tilde{X} , \tilde{Y} , \tilde{F} \right) & = & e^{i} \left( X \right)e^{j} \left( X \right)\tilde{F}^{k}_{i}\tilde{F}^{l}_{j} g_{kl} \left( Y \right), \ \forall \left( \tilde{X} , \tilde{Y} , \tilde{F} \right) .\\
\end{eqnarray*} 
In this way, 

\begin{eqnarray*}
\left. \dfrac{\partial r}{\partial x^{k}} \right | _{j_{X,Y}^{1}\phi} & = & 2 \left. \dfrac{\partial e^{i}}{\partial x^{k} } \right | _{X} e^{j} \left( X \right)\tilde{F}^{k}_{i}\tilde{F}^{l}_{j} g_{kl} \left( Y \right).\\
\end{eqnarray*} 
Hence, by using the non-degeneracy of $g$ we have that $ \dfrac{\partial r}{\partial x^{k}} = 0$ if, and only if,
\begin{equation}\label{106}
\dfrac{\partial e^{i}}{\partial x^{k}} = 0, \ \forall i.
\end{equation}
With this, $\textbf{(1)''}$ is satisfied if, and only if, the vector field $e$ is constant respect to $\left( x^{i} \right)$, i.e.,
\begin{equation}\label{107}
e = \lambda^{i} \dfrac{\partial}{\partial x^{i}}, \ \lambda^{i} \equiv Const.
\end{equation}
Next, we will study condition $\textbf{(2)''}$. Notice that,
$$ \dfrac{\partial J}{\partial x^{i}} = \dfrac{\partial J}{\partial F^{l}_{m}}\dfrac{\partial F^{l}_{m}}{\partial x^{i}}.$$
Using the derivative of $J$ (which we have shown above), we have that
$$\dfrac{\partial J}{\partial {F^{l}_{m}}_{|\tilde{F}}} = det \left( \tilde{F} \right) \left( \tilde{F}^{-1} \right)^{l}_{m}.$$
Then, $\textbf{(2)''}$ is satisfied if, and only if,
\begin{equation}
\dfrac{\partial F^{l}_{m}}{\partial x^{i}} = 0, \ \forall i,l,m.
\end{equation}
Observe that

\begin{eqnarray*}
\left. \dfrac{\partial F^{l}_{m}}{\partial x^{k}} \right | _{j_{X,Y}^{1}\phi} &=& \left. \dfrac{\partial F^{l}_{m}\circ \left( x^{i} , y^{j} , y^{j}_{i} \right)^{-1}}{\partial X^{k}} \right | _{\left( \tilde{X} , \tilde{Y}, \tilde{F}\right)}\\
&=&  \mbox{\scriptsize $\left. \dfrac{\partial}{\partial X^{k}}   \left( \left. \dfrac{\partial X^{l} \circ \left(y^{j}\right)^{-1}}{\partial X^{k}} \right | _{\tilde{Y}} \cdot \tilde{F}^{k}_{r} \cdot \left[ \left. \dfrac{\partial X^{r} \circ \left(x^{i}\right)^{-1}}{\partial X^{m}} \right | _{ \tilde{X} }\right]^{-1}\right)  \right | _{\left( \tilde{X} , \tilde{Y}, \tilde{F}\right)}$}.\\
\end{eqnarray*} 
i.e.,
$$\left. \dfrac{\partial F^{l}_{m}}{\partial x^{k}} \right | _{j_{X,Y}^{1}\phi} = 0,$$
if, and only if,
$$  \left. \dfrac{\partial }{\partial X^{k}}\left( \left. \dfrac{\partial X^{m} \circ \left(x^{i}\right)^{-1}}{\partial X^{i}}\right | _{\tilde{X}}  \right)\right | _{\tilde{X}} = 0.$$
So, $\textbf{(2)''}$ is tantamount to,
$$\dfrac{\partial X^{m}}{\partial x^{i}} \equiv Const, \ \forall i,m.$$
This fact implies that,
$$ e \left( X^{m} \right) \equiv Const, \ \forall m.$$
i.e.,
$$ e = \mu^{i} \dfrac{\partial}{\partial X^{i}}, \ \mu^{i} \equiv Const.$$
Notice that, by using Eq. (\ref{107}), this implies, indeed, that the canonical basis is a (global) system of homogeneous coordinates on $\mathcal{B}$. So, we extract the following conclusions 
\begin{itemize}
\item[\textbf{(a)}] $\mathcal{B}$ is (locally) homogeneous if, and only if, the vector field $e$ is constant respect to the canonical basis of $\mathbb{R}^{3}$. 
\item[\textbf{(b)}] The homogeneity of $\mathcal{B}$ implies that the canonical coordinates are homogeneous coordinates.
\item[\textbf{(c)}] $\mathcal{B}$ is locally homogeneous if, and only if, $\mathcal{B}$ is global homogeneous.

\end{itemize}

\section*{Acknowledgements}
This work has been partially supported by MINECO Grants MTM2016-76-072-P and the ICMAT Severo Ochoa projects SEV-2011-0087 and SEV-2015-0554. V.M.~Jim{\'e}nez wishes to thank MINECO for a FPI-PhD Position. We would like to thank the referees for their valuable suggestions that have contributed to improve this paper.

\bibliographystyle{plain}
\bibliography{Library}

\begin{thebibliography}{10}

\bibitem{DABIL}
{B. A. Bilby}.
\newblock Continuous distributions of dislocations.
\newblock In {\em Progress in solid mechanics, {V}ol. 1}, pages 329--398.
  North-Holland Publishing Co., Amsterdam, 1960.

\bibitem{FBLOO}
F.~{Bloom}.
\newblock {\em Modern differential geometric techniques in the theory of
  continuous distributions of dislocations}, volume 733 of {\em Lecture Notes
  in Mathematics}.
\newblock Springer, Berlin, 1979.

\bibitem{COLE}
B.~D. Coleman.
\newblock Simple liquid crystals.
\newblock {\em Archive for Rational Mechanics and Analysis}, 20:41--58, jan
  1965.

\bibitem{DIGFMCORD}
L.~A. Cordero, C.~T. Dodson, and M.~de~Le{\'o}n.
\newblock {\em Differential Geometry of Frame Bundles}.
\newblock Mathematics and Its Applications. Springer Netherlands, Dordrecht,
  1988.

\bibitem{MELZA}
M.~{El\.zanowski}, M.~Epstein, and J.~{\'Sniatycki}.
\newblock {$G$}-structures and material homogeneity.
\newblock {\em J. Elasticity}, 23(2-3):167--180, 1990.

\bibitem{MELZASEG}
M.~{El\.zanowski} and S.~{Prishepionok}.
\newblock Locally homogeneous configurations of uniform elastic bodies.
\newblock {\em Rep. Math. Phys.}, 31(3):329--340, 1992.

\bibitem{EPSBOOK}
M.~Epstein.
\newblock {\em The Geometrical Language of Continuum Mechanics}.
\newblock Cambridge University Press, 2010.

\bibitem{MEPMDL}
M.~Epstein and M.~de~Le{\'o}n.
\newblock Geometrical theory of uniform {C}osserat media.
\newblock {\em Journal of Geometry and Physics}, 26(1):127 -- 170, 1998.

\bibitem{MEPMDLTHIRD}
M.~Epstein and M.~de~Le{\'o}n.
\newblock Continuous distributions of inhomogeneities in liquid-crystal-like
  bodies.
\newblock {\em Proceedings of the Royal Society of London A: Mathematical,
  Physical and Engineering Sciences}, 457(2014):2507--2520, 2001.

\bibitem{MEPMDLSEG}
M.~Epstein and M.~{de Le{\'o}n}.
\newblock Unified geometric formulation of material uniformity and evolution.
\newblock {\em Math. Mech. Complex Syst.}, 4(1):17--29, 2016.

\bibitem{MGEOEPS}
M.~Epstein, V.~M. Jim{\'e}nez, and M.~de~Le{\'o}n.
\newblock Material geometry.
\newblock {\em Journal of Elasticity}, pages 1 -- 24, Oct 2018.

\bibitem{JDES}
J.~D. Eshelby.
\newblock The force on an elastic singularity.
\newblock {\em Philosophical Transactions of the Royal Society of London A:
  Mathematical, Physical and Engineering Sciences}, 244(877):87--112, 1951.

\bibitem{MD}
V.~M. Jim{\'e}nez, M.~{de Le{\'o}n}, and M.~Epstein.
\newblock Material distributions.
\newblock {\em Mathematics and Mechanics of Solids}, 0(0):1081286517736922, 0.

\bibitem{CHARDIST}
V.~M. Jim{\'e}nez, M.~de~Le{\'o}n, and M.~Epstein.
\newblock Characteristic distribution: An application to material bodies.
\newblock {\em Journal of Geometry and Physics}, 127:19 -- 31, 2018.

\bibitem{COSVME}
V.~M. Jim{\'e}nez, M.~de~Le{\'o}n, and M.~Epstein.
\newblock Lie groupoids and algebroids applied to the study of uniformity and
  homogeneity of cosserat media.
\newblock {\em International Journal of Geometric Methods in Modern Physics},
  15(08):1830003, 2018.

\bibitem{KKON}
K.~Kondo.
\newblock Geometry of elastic deformation and incompatibility.
\newblock 1:5–17, 1955.

\bibitem{EKRONFI}
E.~Kr{\"o}ner.
\newblock {\em Allgemeine Kontinuumstheorie der Versetzungen und
  Eigenspannungen}, volume~4.
\newblock 1960.

\bibitem{EKRON}
E.~{Kr{\"{o}}ner}.
\newblock {\em Mechanics of {Generalized Continua}}.
\newblock Springer, Heidelberg, 1968.

\bibitem{RLAND}
R.~W. Lardner.
\newblock {\em Mathematical Theory of Dislocations and Fracture}.
\newblock Mathematical expositions. University of Toronto Press, Toronto, 1974.

\bibitem{KMG}
K.~C.~H. {Mackenzie}.
\newblock {\em General theory of {L}ie groupoids and {L}ie algebroids}, volume
  213 of {\em London Mathematical Society Lecture Note Series}.
\newblock Cambridge University Press, Cambridge, 2005.

\bibitem{JEMARS}
J.~E. {Marsden} and T.~J.~R. {Hughes}.
\newblock {\em Mathematical foundations of elasticity}.
\newblock Dover Publications, Inc., New York, 1994.
\newblock Corrected reprint of the 1983 original.

\bibitem{GAMAU2}
G.~A. {Maugin}.
\newblock {\em Material inhomogeneities in elasticity}, volume~3 of {\em
  Applied Mathematics and Mathematical Computation}.
\newblock Chapman \& Hall, London, 1993.

\bibitem{IMJMO}
I.~{Moerdijk} and J.~{Mr\v cun}.
\newblock On integrability of infinitesimal actions.
\newblock {\em Amer. J. Math.}, 124(3):567--593, 2002.

\bibitem{FRNAV}
F.~R.~N. Nabarro.
\newblock {\em Theory of crystal dislocations}.
\newblock Dover Books on Physics and Chemistry. Dover Publications, New York,
  1987.

\bibitem{WNOLL}
W.~{Noll}.
\newblock Materially uniform simple bodies with inhomogeneities.
\newblock {\em Arch. Rational Mech. Anal.}, 27:1--32, 1967/1968.

\bibitem{CTRUETOU}
J.~L. Synge.
\newblock {\em Principles of Classical Mechanics and Field Theory}.
\newblock Number v. 3,n.{\textordmasculine} 1 in Handbuch der Physik. Springer,
  Berlin, 1960.

\bibitem{CTRUE}
C.~{Truesdell} and W.~{Noll}.
\newblock {\em The non-linear field theories of mechanics}.
\newblock Springer-Verlag, Berlin, third edition, 2004.
\newblock Edited and with a preface by Stuart S. Antman.

\bibitem{JNM}
J.~N. {Vald{\'e}s}, {\'A}.~F.~T. {Villal{{\'o}}n}, and J.~A.~V.
  {Alarc{{\'o}}n}.
\newblock {\em Elementos de la teor{\'\i}a de grupoides y algebroides}.
\newblock Universidad de C{\'a}diz, Servicio de Publicaciones, C{\'a}diz, 2006.

\bibitem{CCWANTHIRD}
C.~C. Wang.
\newblock A general theory of subfluids.
\newblock {\em Archive for Rational Mechanics and Analysis}, 20:1--40, jan
  1965.

\bibitem{CCWANSEG}
C.~C. {Wang}.
\newblock On the geometric structures of simple bodies. {A} mathematical
  foundation for the theory of continuous distributions of dislocations.
\newblock {\em Arch. Rational Mech. Anal.}, 27:33--94, 1967/1968.

\bibitem{CCWAN}
C.~C. {Wang} and C.~{Truesdell}.
\newblock {\em Introduction to rational elasticity}.
\newblock Noordhoff International Publishing, Leyden, 1973.
\newblock Monographs and Textbooks on Mechanics of Solids and Fluids: Mechanics
  of Continua.

\bibitem{WEINSGROUP}
A.~Weinstein.
\newblock Groupoids: unifying internal and external symmetry. {A} tour through
  some examples.
\newblock In {\em Groupoids in analysis, geometry, and physics ({B}oulder,
  {CO}, 1999)}, volume 282 of {\em Contemp. Math.}, pages 1--19. Amer. Math.
  Soc., Providence, RI, 2001.

\end{thebibliography}
\end{document}